\documentclass[10pt,a4paper,authoryear]{article}
\usepackage{amsfonts,amsthm,amsmath,mathrsfs} 
\usepackage{amssymb,nicefrac,bm}
\usepackage{graphicx}
\usepackage{setspace} 
\usepackage{multirow}
\usepackage{xcolor,hyperref}
\usepackage[utf8]{inputenc}
\usepackage[affil-it]{authblk}
\usepackage[comma]{natbib}

\newtheorem{theorem}{Theorem} %[section]
\newtheorem{proposition}[theorem]{Proposition}
\newtheorem{lemma}[theorem]{Lemma}

\theoremstyle{definition}\newtheorem{remark}[theorem]{Remark}
\theoremstyle{definition}\newtheorem{example}[theorem]{Example}

\newcommand{\bN}{\mathbb{N}}
\newcommand{\bP}{\mathbb{P}}

\newcommand{\bR}{\mathbb{R}}
\newcommand{\bC}{\mathbb{C}}

\newcommand{\bH}{\mathbb{H}}

\newcommand{\bS}{\mathbb{S}}

\newcommand{\cA}{\mathcal{A}}
\newcommand{\cB}{\mathcal{B}}

\newcommand{\cD}{\mathcal{D}}
\newcommand{\cE}{\mathcal{E}}
\newcommand{\cF}{\mathcal{F}}

\newcommand{\cL}{\mathcal{L}}

\newcommand{\unif}{\text{\rm unif}}

\newcommand{\Exp}{\text{\rm Exp}}

\newcommand{\Beta}{\text{\rm Beta}}

\newcommand{\OO}{\mathbb{O}}
\newcommand{\Sd}{\mathbb{S}_{d}}
\newcommand{\Bd}{\mathbb{B}_{d}}

\newcommand{\Seins}{\mathbb{S}_{1}}
\newcommand{\Snull}{\mathbb{S}_{0}}

\renewcommand{\t}{^{\text{t}}}

\newcommand{\SP}{\text{\rm SP}}
\newcommand{\vMF}{\text{\rm MF}}

\newcommand{\CI}{\text{\rm CI}}
\newcommand{\CII}{\text{\rm CII}}
\newcommand{\SN}{\text{\rm SN}}

\newcommand{\Wat}{\text{\rm Wat}}
\newcommand{\SBeta}{\text{\rm SBeta}}
\newcommand{\CHS}{\text{\rm CHS}}
\newcommand{\HS}{\text{\rm HS}}
\newcommand{\chs}{\text{\rm chs}}
\newcommand{\hs}{\text{\rm hs}}
\newcommand{\AG}{\text{\rm AG}}
\newcommand{\DPC}{\text{\rm DPC}}
\newcommand{\dpc}{\text{\rm dpc}}

\newcommand{\geo}{\text{\rm geo}}

\newcommand{\psk}{\text{\rm psk}}
\newcommand{\NB}{\text{\rm NB}}
\newcommand{\nb}{\text{\rm nb}}

\newcommand{\znb}{\text{\rm znb}}

\newcommand{\gpsk}{\text{\rm gpsk}}
\newcommand{\WC}{\text{\rm WC}}
\newcommand{\WN}{\text{\rm WN}}
\newcommand{\br}{\text{\rm br}}

\newcommand{\nf}{\nicefrac}
\allowdisplaybreaks

\begin{document}

\title{Discrete mixture representations of spherical distributions}

\author{Ludwig Baringhaus\thanks{Corresponding author}~~and Rudolf Gr{\"u}bel}

\affil{
  Institute of Actuarial and Financial Mathematics\\
  Leibniz Universit\"at Hannover\\
  Postfach 60 09, D-30060 Hannover, Germany\\
  lbaring@stochastik.uni-hannover.de\\
  rgrubel@stochastik.uni-hannover.de
}

\date{}

\maketitle
\renewcommand{\thefootnote}{\arabic{footnote}}

\begin{abstract}
We obtain discrete mixture representations for parametric families of 
probability distributions on Euclidean 
spheres, such as the von Mises--Fisher, the Watson and the angular Gaussian families. 
In addition to several special results we present a general approach to isotropic 
distribution families that is based on 
density expansions in terms of special surface harmonics. We discuss the connections
to stochastic processes on spheres, in particular random walks, discrete
mixture representations derived from spherical diffusions, 
and the use of Markov representations for the mixing 
base to obtain representations for families of spherical distributions.


\medskip
\textit{MSC2020-Mathematics Subject Classification:} Primary 62H11; secondary 60E05

\end{abstract}

\section{Introduction}\label{sec:intro}
A \emph{discrete mixture representation} for a parametric family $\{P_\theta:\,\theta\in\Theta\}$
of probability measures in terms of another family $\{Q_n:\, n\in \bN_0\}$ of probability
measures, the \emph{mixing base}, all defined on the same measurable space, is of the form
\begin{equation}\label{eq:repr0}
      P_\theta =\sum_{n=0}^\infty w_\theta(n) \, Q_n \quad\text{for all }\theta\in\Theta.
\end{equation}
Here, for each $\theta\in\Theta$, the \emph{mixing coefficients} $(w_\theta(n))_{n\in\bN_0}$
are the individual probabilities of a distribution $W_\theta$, the \emph{mixing distribution},
on (the set of subsets of) $\bN_0$. A classical case is the representation of non-central
chisquared distributions with $k$ degrees of freedom, $P_\theta=\chi_k^2(\theta^2)$ with non-centrality
parameter $\theta^2>0$, as Poisson mixtures of central chisquared distributions, where
$Q_n=\chi_{2n+1}^2:=\chi^2_{2n+1}(0)$ and where $W_\theta$ is the Poisson distribution with mean
$\lambda=\theta^2/2$; see e.g. the books of \citet{LieseMiescke}
and \citet{MoePeres} where the representation appears in statistics in connection with the 
power of statistical tests and in probability theory in connection with the local times 
of Markov processes respectively.
 Such mixture representations can be related to two-stage experiments: 
In order to obtain a value $x$ with distribution $P_\theta$ we first choose $n$ 
according to $W_\theta$ and then choose $x$ according to~$Q_n$. This leads to an
immediate application of discrete mixture representations in the context of 
simulation methodology. 

In the present paper we continue our previous investigations (see \citet{BaGr5, BaGr6}),
and now specifically consider distributions on the Euclidean sphere 
$\bS_d:=\{x\in\bR^{d+1}:\, \|x\|=1\}$ of $(d+1)$-dimensional real vectors of unit 
length. This case seems to us to deserve some interest, 
in particular if specific properties of spheres are taken into account: The group $\OO(d+1)$ 
of orthogonal transformations of the ambient space $\bR^{d+1}$ 
acts transitively on $\bS_d$,
and there is a `polar decomposition' (or `tangent-normal decomposition', see Section \ref{subsec:asrepr})
that relates $\bS_d$ to $[-1,1]\times \bS_{d-1}$
if $d>1$. 

We generally assume that the distribution parameters in the above general setup are of the
form $\theta=(\eta,\rho)$, where $\eta\in\Sd$ may be seen as a location parameter; instead of
$P_\theta$ we also write $P_{\eta,\rho}$. We obtain mixture representations that split the dependence 
on the two parts of the parameter in the sense that
\begin{equation}\label{eq:repr}
      P_{\eta,\rho} =\sum_{n=0}^\infty w_\rho(n) \, Q_{n,\eta}
              \quad\text{for all }\theta=(\eta,\rho)\in\Theta.
\end{equation}
In particular, the mixing distributions depend on $\rho$ only. 
For fixed $\rho$ on the left, or fixed $n\in\bN_0$ on the right hand side of~\eqref{eq:repr},
the families $\{P_{\eta,\rho}:\, \eta\in\bS_d\}$ respectively $\{Q_{n,\eta}:\,\eta\in\bS_d\}$
are parametrized by the sphere and are defined on its Borel subsets $\,\cB(\bS_d)$. 
We assume that these
families interact with the group action mentioned above in the sense that they are 
 \emph{isotropic}; see~\eqref{eq:isotropic} below. 
In particular, their elements are then rotationally symmetric about the axis specified by $\eta$. 
As a simple application of the representation \eqref{eq:repr} we mention that with the finite
sums $R_{\eta,\rho}(A)\,:=\,\sum_{k=0}^n w_\theta(n) \,Q_k(A)$ we have monotonically increasing
approximations of the probabilities $P_{\eta,\rho}(A)$, $A\in \cB(\Sd)$, with uniform error 
bounds in the sense that
\begin{displaymath}
  0\,\le\, P_{\eta,\rho}(A)\,-R_{\eta,\rho}(A)\,\le \sum_{k=n+1}^\infty w_\rho(n)\quad \text{for all}~\eta\in \Sd,\,A\in\cB(\Sd).
\end{displaymath}
The literature contains several other applications; see for example the relation to
nonparametric Bayesian inference in~\citet[Section 5.1]{BaGr6}. 

%The paper is organized as follows.

In Section~\ref{sec:special} we collect some basic notation 
and obtain mixture representations for the von Mises--Fisher family and two spherical
Cauchy families in Theorem~\ref{thm:mixrepr_vMF}, the Watson 
family in Theorem~\ref{thm:mixrepr_Wat}, and an angular Gaussian 
family in Theorem~\ref{thm:mixrepr_Saw}. 
The mixing bases are chosen specifically for the respective 
family, with a view towards reflecting its properties. A different base
will generally lead to a different representation, as demonstrated by~\citet{BaGr5} 
in the context of non-central chisquared distributions. In 
Section~\ref{sec:ultra} we present a general approach that uses
expansions of densities in terms of special surface harmonics. The resulting 
mixing base has a structural property that we call \emph{self-mixing stability}. 
This property makes it comparably easy to relate different expansions to each other. 
We obtain representations with this mixing base for the wrapped Cauchy and the
wrapped normal families in Theorem~\ref{thm:wrapped}, and for the von Mises--Fisher 
families in Theorem~\ref{thm:sphericalMF}. These results only hold under conditions 
on the parameter $\rho$ that ensure that the respective distribution is not too far
away from the uniform distribution on the sphere; in Example~\ref{ex:lincombbase} 
we work out a possibility for extending this range. 
 
For fixed $\rho$ or $n$ we may regard  $P_{\eta,\rho}$ and $Q_{n,\eta}$
as probability kernels via $(\eta,A)\mapsto P_{\eta,\rho}(A)$,
$(\eta,A)\mapsto Q_{n,\eta}(A)$. This 
provides a general connection with Markov processes.
We briefly return to the classical mixture representation of 
non-central one-dimensional chisquared distributions, which may be written as \begin{equation}\label{eq:chiquadrat}
   (X + \theta)^2\; =_{\mathcal D} \, X^2+2\,\sum_{j=1}^{N(\theta)} E_j.
\end{equation}
Here the random variables $N(\theta),X,E_1,E_2,\dots$ are independent, 
$X$ has the standard normal distribution, $E_1,E_2,\dots$ are exponentially distributed 
with mean 1, and $N(\theta)$ has the Poisson distribution with parameter $\theta^2/2$. 
The path-wise point of view displays the distributions $\chi^2_1(\theta)$, $\theta\ge 0$,
as the distributions of randomly stopped partial sums of independent
random variables, and~\eqref{eq:chiquadrat} may be used to read off stochastic 
monotonicity and infinite divisibility of non-central chisquared distributions. 
Note that the representation only covers the one-dimensional marginal distributions
of the process $((X+\theta)^2)_{\theta\ge 0}$, as the left hand side of~\eqref{eq:chiquadrat}
is obviously not pathwise monotone in the `time parameter' $\theta$.
 Quite generally, \eqref{eq:repr0} can be related to randomly stopped stochastic 
processes:  If $X=(X_n)_{n\in\bN_0}$ is such that $X_n$ has distribution $Q_n$ 
for all $n\in\bN_0$ then $X_\tau$ has distribution $P_\theta$ if $\tau$ is independent 
of $X$ and has distribution $W_\theta$. 

In Section~\ref{sec:Markov} we discuss several connections between families of
spherical distributions and stochastic processes on spheres.
We consider random walks on spheres in Section~\ref{subsec:rw}, distribution families 
that arise in connection with diffusion processes in Section~\ref{subsec:diff}, 
and the use of Markov representations of the mixing base in connection with almost 
sure representations for distribution families in Section~\ref{subsec:asrepr}. The 
ultraspherical mixing base from Section~3 will be useful at various stages.

For a single transition kernel
we obtain a family $(X^\eta_n)_{n\in\bN_0}$ of Markov chains indexed by their initial
state $\eta$, meaning that $X^\eta_0=\eta$ with probability 1. Isotropy of the kernel 
then extends to isotropy of the corresponding distributions on the path space. 
For the elements of the mixing base in Section~\ref{sec:ultra} the marginal 
distributions of these chains have a particular simple description.  
Further, isotropy relates a family
$\{P_\eta:\,\eta\in\bS_d\}$ to a single distribution on $[-1,1]$ via the 
latitude projection $x\to\eta^t x$, with $\eta$ as `north pole'. For the chain
$X^\eta=(X^\eta_n)_{n\in\bN_0}$ we obtain an associated latitude process $Y=(Y_n)_{n\in\bN_0}$
via $Y_n:=\eta\t X_n$ for all $n\in\bN_0$. For isotropic kernels this is again a Markov
chain, now on $[-1,1]$ and with start at 1.  
Finally, for the von Mises--Fisher distributions we show that a homogeneous 
Markov process on the sphere with these as marginal distributions does not exist, 
see Theorem~\ref{thm:vMF_Markov}, and  we obtain a result similar to~\eqref{eq:chiquadrat}, see Example~\ref{ex:vMFasrepr}.

Proofs are collected in Section~\ref{sec:proofs}.

Mixing of distributions is a standard topic in probability theory and statistics, 
see e.g.~\citet{Lindsay}. Spherical data and families of spherical distributions 
have similarly been investigated for a long time and by many researchers; standard
references are the classic monograph of~\citet{Watson} and, more recently, the book of~\citet{MardiaJupp}.
For a review of distributions on spheres we refer to~\citet{Pewsey}, 
see also~\citet{WatsonJAP}. Of particular interest for the topics treated here
is the very recent paper of~\citet{Mija} where a discrete
mixture representation for the marginal distributions of spherical Brownian motion
is developed. More specific references will be given at the appropriate places below.

%%%%%%%%%%%%%%%%%%%

\section{Generalities and some special results}\label{sec:special}\label{sec:general}
We need some basic notions and definitions. 
We write $X\sim\mu$ if $X$ is a random variable on some background probability  
space $(\Omega,\cF,\bP)$ with distribution $\mu$. Formally, let $(\Omega,\cA)$ and 
$(\Omega',\cA')$ be measurable spaces and suppose that $T:\Omega\to\Omega'$ is 
$(\cA,\cA')$-measurable. Then the push-forward $P^T$ of a probability measure $P$ on 
$(\Omega,\cA)$ under $T$ is the probability measure on $(\Omega',\cA')$
given by $P^T(A)=P(T^{-1}(A))$, $A\in\cA'$, and $X\sim\mu$ is the same as $\bP^X=\mu$.
For many of the measurable spaces considered below there is a canonical
uniform distribution, often defined by invariance under a group operation. 
To avoid tiresome repetitions we agree that densities refer 
to the respective uniform distribution if not specified otherwise. 

We fix a dimension $d\ge 1$, but instead of $d$ we often use 
\begin{equation}\label{eq:lambda}
\lambda=\lambda(d):=(d-1)/2,
\end{equation}
as this is common in connection with families of special functions.
In particular, whenever $d$ and $\lambda$ appear together, they are related by~\eqref{eq:lambda}. 
The group $\OO(d+1)$ of orthogonal $(d+1)\times(d+1)$-matrices $U$ acts on $\Sd$ via
$x\mapsto Ux$, and the uniform distribution $\unif(\Sd)$ on the sphere is the unique
probability measure on the Borel subsets of $\bS_d$ that is invariant under all such
transformations. For a fixed $\eta\in\Sd$ the push-forward $\nu_d$ of $\unif(\Sd)$ 
under the mapping $x\mapsto \eta\t x$ has density $h_d$ with respect to 
the uniform distribution $\unif(-1,1)$ on the interval $[-1,1]$, where
\begin{equation}\label{eq:dens_hd}
h_d(y) := \frac{\Gamma(\lambda+1)}{\Gamma(\nf{1}{2})\Gamma(\lambda+\nf{1}{2})}
                    \,(1-y^2)^{\lambda-1/2}, \quad -1<y<1.
\end{equation}
Note that this does not depend on $\eta\in\Sd$.  
Further, if $X\sim\unif(\Sd)$ and $Y=\eta\t X\sim \nu_d$, then the conditional distribution 
of $X$ given $Y=y$ is the uniform distribution on 
\begin{equation}\label{eq:defC}
      C_d(\eta,y):= \{x\in\Sd:\, \eta\t x=y\},
\end{equation} 
with $\unif(C_d(\eta,y))$ the unique probability measure on this set that is invariant 
under the subgroup $\{U\in\OO(d+1):\, U\eta=\eta\}$ of $\OO(d+1)$. This may be seen in the
context of the polar decomposition mentioned in the introduction.

Conversely, given a probability measure $\nu$ on $[-1,1]$ and a parameter $\eta\in\Sd$, we can
construct a distribution $\mu=\mu_\eta$ on $\Sd$ via the kernel 
$(y,A)\mapsto \unif(C_d(\eta,y))(A)$. In particular, for bounded and measurable
functions $\phi:\Sd\to\bR$,
\begin{equation*}
\int \phi(x)\, \mu(dx) \; 
         =\; \int_{[-1,1]} \int_{C_d(\eta,y)} \phi(x)\, \unif(C_d(\eta,y))(dx)\, \nu(dy).
\end{equation*}  
For $\eta\in\bS_d$ and  a measurable function $g:[-1,1]\to \bR$
the function $f_\eta:\Sd \to \bR$ given by 
\begin{equation}\label{eq:basic}
  f_\eta(x)=g(\eta\t x),\quad x\in\Sd,
\end{equation}
is $\unif(\bS_d)$-integrable if and only if $g$ is $\nu_d$-integrable, and then
\begin{equation*}
  \int f_\eta(x)\,\unif(\bS_d)(dx) = \int g(t)\, \nu_d(dt).
\end{equation*}
Thus $f_\eta$ is the density of a probability measure on $\Sd$ if and only if
$g$ is the $\nu_d$-density of a probability measure on $[-1,1]$. In particular,
a probability density $g$ on $[-1,1]$ generates a family $\{Q_\eta:\, \eta\in\bS_d\}$ 
of spherical distributions via~\eqref{eq:basic}, and such families are 
\emph{isotropic} in the sense that 
\begin{equation}\label{eq:isotropic}
    Q_{\eta}^U=Q_{U\eta}\quad\text{for all } U\in \OO(d+1),\,\eta\in\bS_d. 
\end{equation}
In particular, each $Q_\eta$ is invariant under all rotations with axis~$\eta$. As the function
$\eta\mapsto \int_A g(\eta\t x)\,\unif(\Sd)(dx)$ is $\cB(\Sd)$-measurable for all $A\in\cB(\Sd)$,
$Q_{\bm\cdot}:\Sd\times \cB(\Sd)\rightarrow \bR$ defined by
$Q_{\bm\cdot}(\eta,A)=Q_\eta(A),\,(\eta,A)\in \Sd\times \cB(\Sd),$ is a Markov kernel from
$\left (\Sd,\cB(\Sd)\right )$ to $\left (\Sd,\cB(\Sd)\right )$. Further, if $Q$ is a Markov kernel from
$\left (\Sd,\cB(\Sd)\right )$ to $\left (\Sd,\cB(\Sd)\right )$ and $U\in \OO(d+1)$,
then the kernel $Q^U:\Sd\times \cB(\Sd)\rightarrow \bR$ defined by
$Q^U(\eta,A)=Q(\eta,U\t A),\,(\eta,A)\in \Sd\times \cB(\Sd),$ with
$U\t A:=\{U\t x:x\in A\}$, is the push-forward of $Q$ under $U$. The kernel $Q$ is 
isotropic if
\begin{equation}\label{ker:isotropic}
    Q^U(\eta,\bm \cdot )=Q(U\eta,\bm \cdot )\quad\text{for all}~\eta\in \Sd~\text{and all}~U\in \OO(d+1).
\end{equation}
Some classical special functions will be needed below. Let
\begin{equation*}
(\alpha)_n:=\frac{\Gamma(\alpha+n)}{\Gamma(\alpha)},
    \quad\alpha>0,\,n\in\bN_0,
\end{equation*} 
be the ascending factorials. The \emph{modified Bessel functions $I_\alpha$ of the first kind}
are given by
\begin{equation}\label{BesselI}
    I_\alpha(x) = \sum_{n=0}^\infty \frac{1}{n! \,
                      \Gamma(n+\alpha+1)} \, \Bigl(\frac{x}{2}\Bigr)^{2n+\alpha}, 
                      \quad x\ge 0,
\end{equation}        
with real nonnegative parameter $\alpha$,
the \emph{confluent hypergeometric functions} are
\begin{equation*}
   {}_1F_1(\alpha;\beta;x) = \sum_{n=0}^\infty \frac{(\alpha)_n}{(\beta)_n\, n!} \, x^n,\quad x\in\bR,      
\end{equation*}
with real positive parameters $\alpha$ and $\beta$, and the \emph{hypergeometric functions} are
\begin{equation*}
   {}_2F_1(\alpha,\beta;\gamma;x) = \sum_{n=0}^\infty \frac{(\alpha)_n(\beta)_n}{(\gamma)_n\,n!} \, x^n,\quad -1<x<1,      
\end{equation*}
with real positive parameters $\alpha$, $\beta$, and $\gamma$.

\begin{example}\label{ex:bsp} 
(a) Let $p>-\nf{1}{2}$ and let 
$\nu$ be the distribution on $[-1,1]$ with  $\unif(-1,1)$-density 
$y\mapsto \frac{\Gamma(p+\lambda+1)}{\Gamma(p+\nf{1}{2})\Gamma(\lambda+\nf{1}{2})}|y|^{2p}(1-y^2)^{\lambda-1/2}$.
Then $\nu$ has $\nu_d$-density $y\mapsto \frac{\Gamma(\nf{1}{2})\Gamma(p+\lambda+1)}
{\Gamma(p+\nf{1}{2})\Gamma(\lambda+1)}\,|y|^{2p},$ and we obtain the
\emph{spherical power} distribution $\SP_d(\eta,p)$ with density
\begin{equation*}
    f_d^{\SP}(x|\eta,p)\, =\, \frac{\Gamma(\nf{1}{2})\Gamma(p+\lambda+1)}
     {\Gamma(p+\nf{1}{2})\Gamma(\lambda+1)}\,|\eta\t x|^{2p}, \quad x\in\Sd.
\end{equation*}
We mainly use this with $p=n\in\bN_0$, and then have
\begin{equation*}
    f_d^{\SP}(x|\eta,n)\, =\, \frac{(\lambda+1)_n} {(\nf{1}{2})_n}\,
                      (\eta\t x)^{2n}, \quad x\in\Sd.
\end{equation*}

\vspace{.5mm}
(b) Starting with $\nu=\Beta_{[-1,1]}(p+\lambda-\nf{1}{2},q+\lambda-\nf{1}{2})$, $p,q>\nf{1}{2}-\lambda$, the beta distributions on $[-1,1]$
with densities $y\mapsto c(p+\lambda-\nf{1}{2},q+\lambda-\nf{1}{2}) (1-y)^{p+\lambda-\nf{3}{2}}(1+y)^{q+\lambda-\nf{3}{2}}$, where
$c(p+\lambda-\nf{1}{2},q+\lambda-\nf{1}{2})=\Gamma(p+q+2\lambda-1)/\left (2^{p+q+2(\lambda-1)}\Gamma(p+\lambda-\nf{1}{2})\Gamma(q+\lambda-\nf{1}{2})\right )$,
we obtain the \emph{spherical beta} distributions $\SBeta_d(\eta,p,q)$, with densities
\begin{equation*}
f_d^\SBeta(x|\eta,p,q)\, = \, c_d(p,q)\, (1-\eta\t x)^{p-1}(1+\eta\t x)^{q-1},\quad
    x\in\Sd,
\end{equation*}
where the norming constants are given by
\begin{equation}\label{eq:normSBeta}
      c_d(p,q)\; =\; 2^{-\left (p+q+2(\lambda-1)\right )}
                  \frac{\Gamma(\nf{1}{2})\Gamma(\lambda+\nf{1}{2})\Gamma(p+q+2\lambda-1)}
                    {\Gamma(\lambda+1)\Gamma(p+\lambda-\nf{1}{2})\Gamma(q+\lambda-\nf{1}{2})}.
\end{equation} 
%For $U\in\OO(d)$ it holds that $\SBeta_d(\eta,p,q)^U=\SBeta_d(U\eta,p,q)$, hence the
%norming constant does indeed not depend on $\eta$.

\vspace{.5mm}
(c) The \emph{von Mises--Fisher} distributions, which we denote by
$\vMF_d(\eta,\rho)$, $\rho > 0$, arise if we start with $\nu_d$-density proportional 
to $y\mapsto\exp(\rho y)$, $-1\le y\le 1$. 
The continuous density of the associated spherical distribution is
\begin{equation*}
f_d^\vMF(x|\eta,\rho)\, = \, c_d(\rho)\, \exp(\rho\, \eta\t x), \quad x\in\Sd,
\end{equation*}
where the norming constants are given by
\begin{equation*}
         c_d(\rho)\; =\; \frac{\rho^\lambda}
         {2^\lambda\Gamma(\lambda+1)I_\lambda(\rho)}\,.
\end{equation*} 
Further, $\vMF_d(\eta,0)=\unif(\Sd)$. 

\vspace{.5mm}
(d)  The \emph{Watson} distributions
$\Wat_d(\eta,\rho)$, $\rho\in\bR$, arise if we begin with $\nu_d$-density proportional 
to $y\mapsto\exp(\rho y^2)$, $-1\le y\le 1$. The continuous density of the associated
spherical distribution is
\begin{equation*}
f^\Wat_d(x|\eta,\rho)\, = \, c_d(\rho)\, \exp\bigl(\rho\, (\eta\t x)^2\bigr), 
            \quad x\in\bS_d,
\end{equation*}
with norming constants 
$\ c_d(\rho) =  \bigl({}_1F_1(\nf{1}{2};\lambda+1;\rho)\bigr)^{-1}$. 
Clearly, $\Wat_d(\eta,0)=\unif(\Sd)$. 

\vspace{.5mm}
(e) The \emph{angular Gaussian} distributions are the distributions of $X=Z/\| Z \|$, where
$Z$ has the $(d+1)$-variate normal distribution $N_{d+1}(a,\Sigma)$ with mean vector
$a\in\bR^{d+1}\setminus\{0\}$ and symmetric positive definite covariance matrix $\Sigma$; 
see, e.g. \citet[p.\,108]{Watson}. Here, we exclusively deal with the case where $\Sigma$ is
the identity matrix $I_{d+1}$, as the radial parts then lead to isotropic families. 
The distributions arising in this special case seem to have first been studied
in detail by \citet{Saw}. Putting $\eta = a/\|a\|$ and
$\rho = (\frac{1}{2}\|a\|^2)^{\nf{1}{2}}$ we denote by $\AG_d(\eta,\rho)$ the
distribution of $X$ and speak of the angular Gaussian distribution with parameters $\eta$ and
$\rho$. Its density is represented by the infinite series
\begin{equation*}
  f^\AG_d(x|\eta,\rho)\, =e^{-\rho^2}\sum_{k=0}^\infty (2\rho\,\eta\t x)^k \frac{\Gamma((d+1+k)/2)}{k!\Gamma((d+1)/2)}, \quad x\in \bS_d,\, \rho>0.
\end{equation*}
We refer to \citet{Saw}, where it is also pointed out that, with a random variable 
$S\sim \chi_{d+1}^2$, the density can be written as 
\begin{equation}\label{Saw:rep0}
  f^\AG_d(x|\eta,\rho)\, = E\left (e^{-\rho^2}\exp( 2^{\nf{1}{2}}\rho\,S^{\nf{1}{2}}\,\eta\t x)\right ).
\end{equation}

\vspace{.5mm}
(f) We consider two types of spherical Cauchy distributions: The 
\emph{spherical Cauchy distributions of type I} have the $\unif(\Sd)$-densities
\begin{equation*}
  f_d^{\CI}(x|\eta,\rho)\,=\left (\frac{1-\rho^2}{1-2\rho\,\eta\t\,x +\rho^2}\right )^d,\quad x\in\bS_d,
\end{equation*}
and the \emph{spherical Cauchy distributions  of type II} have the $\unif(\Sd)$-densities
\begin{equation*}
  f_d^{\CII}(x|\eta,\rho)\,=\frac{1-\rho^2}{\left (1-2\rho\,\eta\t\,x +\rho^2\right )^{(d+1)/2}},\quad x\in\bS_d,
\end{equation*}
both with parameters $\eta\in \bS_d$ and $\rho\in (0,1)$. We denote by ${\CI}_d(\eta,\rho)$ the
distribution with $\unif(\Sd)$-density  $f_d^{\rm CI}(\,\cdot\,|\eta,\rho)$,
and by ${\CII}_d(\eta,\rho)$ the distribution with
$\unif(\Sd)$-density $f_d^{\rm CII}(\,\cdot\,|\eta,\rho)$. Clearly, ${\CI}_d(\eta,0)\,=\,\CII_d(\eta,0)\,=\,\unif(\Sd)$. 
\end{example}

\vspace{1mm}
The distributions in Example \ref{ex:bsp}~(a) - (e) and their push-forwards under $x\mapsto \eta\t x$, respectively,
are all classical; for basic as well as specific properties
and interesting historical comments we refer to~\citet{WatsonJAP,Watson} and \citet{MardiaJupp}.
It is well known, for example, that the von Mises--Fisher
family in part~(c) arises from the multivariate normal distributions in part~(e)
by conditioning on $\|Z\|$. A well known relation with Brownian motion on $\Sd$
is addressed in Section \ref{subsec:diff} below. In the special case $d=1=(d+1)/2$ the distributions
$\WC_1(\eta,\rho):=\CI_1(\eta,\rho)=\CII_1(\eta,\rho)$
in Example \ref{ex:bsp} (f) are known as the wrapped Cauchy or circular Cauchy distributions;
see, e.g.~\citet{MardiaJupp} and Section \ref{sec:ultra} below.
Hence, with the two types of spherical Cauchy distributions given above we
have two different extensions of this distribution family to higher dimensions. 
For distinction, we
added the name supplement `of type I' and `of type II', respectively. 
The spherical Cauchy distributions of type I were introduced and studied by \citet{Kato}.
Generalizing results obtained by \citet{McCullagh96} for $d=1$, the 
authors especially deal with the behavior of the spherical Cauchy distributions of type I
under M{\"o}bius transformations. The densities $f_d^{\CII}(\,\cdot\,|\eta,\rho)$ were 
considered by
\citet{McCullagh89}, though the author does not speak of spherical Cauchy distributions but,
with the push-forward of $\CII(\eta,\rho)$ under $x\mapsto \eta\t x$, of a noncentral version
of the univariate symmetric beta distribution.

We recall from the introductory remarks that for a discrete mixture representation 
we need a mixing base $(Q_{n,\eta})_{n\in\bN_0}$, $\eta\in\bS_d$, where
each $Q_{n,\eta}$ is a probability measure on the sphere, and mixing distributions
on~$\bN_0$ that depend on $\rho$ only; see~\eqref{eq:repr}. 
Of special interest in the latter context are the
the \emph{negative binomial distributions} $\NB(r,p)$ with parameters $r>0$, $p\in (0,1)$, and
probability mass function
\begin{align*}
      \nb(n|r,p)=\frac{(r)_n}{n!}(1-p)^np^r,\quad n\in\bN_0,
\end{align*}
the 
\emph{confluent hypergeometric series distributions} $\CHS(\alpha,\beta,\tau)$ 
on $\bN_0$ with parameters $\alpha,\beta,\tau>0$ and probability mass functions
\begin{equation*}
 \chs(n|\alpha,\beta,\tau) \,=\, \bigl({}_{1}F_1(\alpha;\beta;\tau)\bigr)^{-1}\;
             \frac{(\alpha)_n}{(\beta)_n}\frac{\tau^n}{n!}, \quad  n\in\bN_0,
\end{equation*}
and the 
\emph{hypergeometric series distributions} $\HS(\alpha,\beta,\gamma,\tau)$ 
on $\bN_0$ with parameters $\alpha,\beta,\gamma,\tau>0$ and probability mass functions
\begin{equation*}
 \hs(n|\alpha,\beta,\gamma,\tau) \,=\, \bigl({}_{2}F_1(\alpha,\beta;\gamma;\tau)\bigr)^{-1}\;
             \frac{(\alpha)_n(\beta)_n}{(\gamma)_n}\frac{\tau^n}{n!}, \quad  n\in\bN_0.
\end{equation*}
For $\tau=0$ we take $\CHS(\alpha,\beta,\tau)$ and $\HS(\alpha,\beta,\tau)$ to be
the one-point mass at 0. The distribution $\CHS(\alpha,\beta,\tau)$
arises as the stationary distribution of a birth-death process with birth rates
$(\alpha + i)\tau$ and death rates $i(\beta + i-1)$, $i\in\bN_0$; 
see \citet{Hall}. Note that these three distribution families are subclasses of the
family of generalized hypergeometric distributions considered recently by
\citet{Themangani}. We also require that each
family $\{Q_{n,\eta}:\, \eta\in\bS_d\}$ is isotropic. 

We can now state our first results. Let $d\in\bN$ be fixed and let $\lambda$ 
be as in~\eqref{eq:lambda}.

\begin{theorem}\label{thm:mixrepr_vMF}
\emph{(a)} The family $\,\{\vMF_d(\eta,\rho):\, \eta\in\Sd, \rho\ge 0\}\,$ has a unique discrete mixture representation 
with mixing base $\SBeta_d(\eta,1,n+1)$, $n\in\bN_0$. This
representation is given by
\begin{equation}\label{eq:repr_vMF}
   \vMF_d(\eta,\rho)\, = \, \sum_{n=0}^\infty \chs(n|\lambda+\nf{1}{2},2\lambda+1,2\rho) \, \SBeta_d(\eta,1,n+1). 
 \end{equation}

\vspace{.5mm}
\noindent
\emph{(b)} The family $\,\{\CI_d(\eta,\rho):\, \eta\in\Sd, \rho \in (0,1)\}\,$  
has a unique discrete mixture representation with mixing base $\SBeta_d(\eta,1,n+1)$, $n\in\bN_0$. This
representation is given by
\begin{equation}\label{eq:repr_CI}
   \CI_d(\eta,\rho)\, = \, \sum_{n=0}^\infty \nb\left (n|\lambda+\nf{1}{2},4\rho/(1+\rho)^2\right ) \, \SBeta_d(\eta,1,n+1). 
 \end{equation}

\vspace{.5mm}
\noindent
\emph{(c)} The family $\,\{\CII_d(\eta,\rho):\, \eta\in\Sd, \rho \in (0,1)\}\,$  
has a unique discrete mixture representation with mixing base $\SBeta_d(\eta,1,n+1)$, $n\in\bN_0$. This
representation is given by
\begin{equation}\label{eq:repr_CII}
   \CII_d(\eta,\rho)\, = \, \sum_{n=0}^\infty \hs\left (n|\lambda+\nf{1}{2},\lambda+1,2\lambda+1,4\rho/(1+\rho)^2\right ) \, \SBeta_d(\eta,1,n+1). 
 \end{equation}
\end{theorem}

In our next result the mixing base depends on the value of $\rho$.

\begin{theorem}\label{thm:mixrepr_Wat}
\emph{(a)} The family $\{\Wat_d(\eta,\rho):\, \eta\in\Sd, \rho \ge 0\}$ 
has a unique discrete mixture representation with mixing base $\SP_d(\eta,n)$, $n\in\bN_0$. This
representation is given by
\begin{equation}\label{eq:repr_Wat1}
   \Wat_d(\eta,\rho)\, = \, \sum_{n=0}^\infty \chs(n|1/2,\lambda+1,\rho) \, \SP_d(\eta,n). 
\end{equation}

\vspace{.5mm}
\noindent
\emph{(b)} The family $\{\Wat(\eta,\rho):\, \eta\in\Sd, \rho \le 0\}$
has a unique discrete mixture representation with mixing base $\SBeta_d(\eta,n+1,n+1)$, $n\in\bN_0$.
This
representation is given by
\begin{equation}\label{eq:repr_Wat2}
   \Wat_d(\eta,\rho)\, = \, \sum_{n=0}^\infty \chs(n|\lambda+1/2,\lambda+1,-\rho) \, \SBeta_d(\eta,n+1,n+1). 
 \end{equation}
\end{theorem}

In order to obtain a similar representation for the family of 
$\{{\AG}_d(\eta,\rho):\, \eta\in\Sd,\rho>0\}$ of angular Gaussian distributions 
we make use of the integral representation 
\begin{equation}\label{def:pc}
  D_\nu(z)=\frac{e^{-z^2/4}}{\Gamma(-\nu)}\int_0^\infty t^{-\nu -1} 
                e^{-zt -t^2/2}\,dt, \quad z\in \bR,
\end{equation}
of the parabolic cylinder functions $D_\nu$ with real index $\nu<0$; 
see \citet[p.\,328]{MOS}.

\begin{lemma}\label{lem:dpc}
Let $\delta>\nf{1}{2}$ and $\tau>0$. Then $\dpc(\,\cdot\,|\delta,\tau)$ with
\begin{equation}\label{eq:dpc}
  \dpc(k|\delta,\tau) \, 
            = \, \frac{1}{k!}(2^{\nf{1}{2}}\tau)^k 2^{k+\delta -1}
                 \frac{(k+2\delta-1)\Gamma(k+\delta-\nf{1}{2})}{\Gamma(\nf{1}{2})}
                           \, e^{-\tau^2/2} D_{-(k+2\delta)}(2^{\nf{1}{2}}\tau),
\end{equation}
$k\in\bN_0$, is a probability mass function.
\end{lemma}

We write $\DPC(\delta,\tau)$ for the associated \emph{discrete parabolic cylinder
distribution} with parameters $\delta>\nf{1}{2}$ and $\tau>0$. For the special 
values $\delta=\lambda + 1=(d+1)/2$ with $d\in\bN$ the statement of the lemma also 
follows from~\eqref{eq:repr_Saw}  below as the values on the right hand side of~\eqref{eq:dpc}
are all nonnegative.

\begin{theorem}\label{thm:mixrepr_Saw}
The family $\{\AG_d(\eta,\,\rho):\, \eta\in\Sd, \rho > 0\}$ 
has a unique discrete mixture representation with mixing base 
$\SBeta_d(\eta,1,n+1)$, $n\in\bN_0$.
This representation is given by
\begin{equation}\label{eq:repr_Saw}
  \AG_d(\eta,\rho)\, = \, \sum_{n=0}^\infty \dpc(n|\lambda+1,\rho) \, \SBeta_d(\eta,1,n+1).
\end{equation}
\end{theorem}

\begin{remark}\label{rem:bem}
(a) Regarding probability measures as real functions on a set of events, we may define the series
in~\eqref{eq:repr_vMF} - \eqref{eq:repr_Wat2} and \eqref{eq:repr_Saw}
as referring to pointwise convergence of functions.
In fact, as the distributions involved all have smooth densities and compact domain, convergence
even holds with respect to uniform convergence in spaces of continuous functions.

(b) The representations are minimal in the sense that the respective mixing base cannot be reduced.
This follows from the uniqueness and the fact that the mixing probabilities are strictly positive.

(c) In connection with the base in Theorem~\ref{thm:mixrepr_vMF} 
all mixtures have densities that are increasing in $\eta\t x$, and in Theorem~\ref{thm:mixrepr_Wat}
and Theorem~\ref{thm:mixrepr_Saw} all mixtures are invariant under the reflection $x\mapsto -x$.     

(d) Interestingly, the von Mises--Fisher family, the spherical Cauchy families and the angular
Gaussian family have the same mixing base of spherical beta distributions. So, these families are obtained
by picking at random (the index $n$ of) the element $\SBeta_d(\eta,1,n+1)$ according to the
respective mixing distributions. Another family with this mixing base is the family of spherical
normal distributions; see Section \ref{subsec:diff}.

(e) For a discussion of other similarities
as well as differences between the von Mises--Fisher family and the spherical Cauchy family
of type I we refer to \citet{Kato}. The von Mises--Fisher family and the spherical Cauchy family
of type II both have representations in terms of multivariate Brownian motion. To be specific,  
let $X=(X_t)_{t\ge 0}$ be a standard Brownian motion in $\bR^{d+1}$, let $Y=(Y_t)_{t\ge 0}$ with
$Y_t=\rho \eta t + X_t$ for $t\ge 0$ be the drifted standard Brownian motion with 
constant drift vector $\rho \eta$, where $\rho\ge 0$, $\eta\in \Sd$, and 
let $Z=(Z_t)_{t\ge 0}$ with $Z_t=\rho\eta +X_t$ for $t\ge 0$, where $0\le \rho<1$, $\eta\in \Sd$, be the Brownian motion starting at $\rho \eta$. With $T_Y:=\inf\{t\ge 0: \|Y_t\|\ge 1\}$
as the first time that $Y$ exits the
Euclidean unit ball $\Bd=\{x\in \bR^{d+1}: \|x\|<1\}$ it then holds that 
$Y_{T_Y}\sim \vMF_d(\eta,\rho)$; see \citet{Gatto} for a more recent proof and 
historical remarks on this result.
Further, with $T_Z:=\inf\{t\ge 0: \|Z_t\|\ge 1\}$ the first time that $Z$ exits 
$\Bd$ it holds that $Z_{T_Z}\sim \CII_d(\eta,\rho)$;
see \citet[p.\,170]{Chung} and \citet{McCullagh89}. 
\end{remark}

\section{Ultraspherical mixing bases}\label{sec:ultra}
Our aim in this section is a mixing base that is applicable 
for general distribution families where, as before, we consider distributions 
$P_\theta$ on $(\bS_d,\cB(\bS_d))$, with $\theta=(\eta,\rho)\in\bS_d\times I$
and $I\subset \bR_+$ an interval,  that have densities $f_\theta$
of the form
\begin{equation}\label{eq:struct}
   f_\theta(x)=g_\rho(\eta\t x),\quad x\in\bS_d.
\end{equation}
We will occasionally omit $d$ or $\lambda$ from the notation.
Recall that $\lambda=(d-1)/2$ and that $\nu_d$ is the push-forward of $\unif(\bS_d)$
under the mapping $x\mapsto \eta\t x$. 
%We also write $\mu_\lambda$ instead of $\nu_d$.OD

We assume that the functions $g_\rho$ in~\eqref{eq:struct}
are elements of
\begin{equation*}
      \bH_\lambda := L^2\bigl([-1,1],\cB([-1,1]),\nu_d\bigr),
\end{equation*}
and on $\bH_\lambda$ we use the inner product
\begin{equation*}
      \langle f,g\rangle_\lambda = \int f(t)g(t)\,\nu_d(dt)
\end{equation*}
and the norm $\|f\|_\lambda=\langle f,f\rangle_\lambda^{1/2}$. 
Then $(\bH_\lambda,\langle \cdot,\cdot\rangle_\lambda)$ is a Hilbert space.
We deal with a special complete sequence of orthogonal polynomials in this
space. For $d=1$ and $\lambda=0$ this is the sequence of
Chebyshev polynomials $T_n$ of the first kind of degree $n\in\bN_0$, 
for $d>1$ and $\lambda >0$ we use the sequence of Gegenbauer or ultraspherical polynomials
$C_n^\lambda$ of degree $n\in\bN_0$; see \citet[Chs. X, XI]{ErdelyiII}. These 
functions play an important role in directional statistics, especially nonparametric directional
statistics; see e.g.\ the papers of \citet{Bingham}, \citet{Gine}, \citet{Prentice}, \citet{Baringhaus},
\citet{Jupp}, and \citet{GPV2021}. 
The functions are standardized such that
\begin{equation*}
  C_n^\lambda(1)=\frac{\Gamma(2\lambda+n)}{\Gamma(2\lambda)n!}
                =\frac{(2\lambda)_n}{n!}\quad \text{for all } n\in\bN_0
\end{equation*}
if $\lambda > 0$; further, $T_n(1)=1$ for all $n\in \bN_0$. 
In particular, $C_0^\lambda\equiv 1 \equiv T_0$. Of course, for $n>0$ none of these functions
is a probability density with respect to $\nu_d$. However, it is known that the Chebyshev polynomials and, for
$\lambda>0$, the Gegenbauer polynomials attain their absolute maximum on $[-1,1]$ at $t=1$; see \citet[p. 206, formula (7)]{ErdelyiII}
and \citet[p.\,786]{AbramSteg}. Hence the standardization 
\begin{equation*}
  D_n^\lambda(t)\, :=\;
  \begin{cases} T_n(t) & \text{if}~\lambda = 0,\\
    C_n^\lambda(1)^{-1}\, C_n^\lambda(t)& \text{if}~\lambda > 0,
  \end{cases}  
\end{equation*}
provides a sequence $(D_n^\lambda)_{n\in\bN_0}$ of orthogonal polynomials that are bounded in absolute
value on $[-1,+1]$ by their value~$1$ in $t=1$. As $(D_n^\lambda)_{n\in\bN_0}$ is complete in~$\bH_\lambda$, we have the
series expansion converging in $\bH_\lambda$ 
\begin{equation}\label{eq:expansion}
      g_\rho \;=\;\sum_{n=0}^\infty \langle g_\rho,D_n^\lambda \rangle_\lambda \,  \|D_n\|_\lambda^{-2}\, D_n^\lambda
     \; =\;  1 +\sum_{n=1}^\infty \beta_n(\rho)\, D_n^\lambda, 
\end{equation}
with $\beta_n(\rho):= \langle g_\rho,D_n^\lambda \rangle_\lambda \,\|D_n\|_\lambda^{-2}$ for
$n\in\bN$. 
We assume throughout this section that $g_\rho$ is such that
\begin{equation}\label{eq:def_beta}
     \beta(\rho):=\sum_{n=1}^\infty |\beta_n(\rho)| \,<\, \infty.
\end{equation}
For fixed $n\in\bN_0$ and $\eta\in\Sd$ the function $H_{n,\eta}^\lambda:\Sd\rightarrow \bR$ defined by
\begin{align*}
  H_{n,\eta}^\lambda(x)=D_n^\lambda(\eta\t x),\,x\in \Sd,
\end{align*}
is the unique \emph{surface harmonic of degree $n$} that depends only on $\eta\t x$ and that
satisfies $H_{n,\eta}^\lambda(\eta)=1$; see \citet[p.\,238]{ErdelyiII}.
Obviously, $H_{n,\eta}^\lambda$
is invariant under the subgroup $\{U\in\OO(d+1):\, U\eta=\eta\}$ of $\OO(d+1)$, 
i.e. for all elements $U$ of the subgroup it holds that $H_{n,\eta}^\lambda(Ux)=H_{n,\eta}^\lambda(x)$ for all
$x\in\Sd$.

Let
\begin{gather}
\begin{aligned}\label{eq:defgamma}
  \gamma^0_0&\,:=1,~\gamma^0_n\,:=1/2\;\quad\text{for all}~n\in\bN,\\
  \gamma^\lambda_n&\, := \, \frac{1}{(1+n/\lambda)\,C_n^\lambda(1)}\quad\text{for all }\lambda>0,~n\in\bN_0.
\end{aligned}
\end{gather}                    
We will repeatedly make use of the basic formulas
\begin{align}
  \int D_m^\lambda(\eta\t x)D_n^\lambda(\xi\t x)\, \unif(\bS_d)(dx)& \ = \ 0,
       \quad \eta,\xi \in \Sd,\,m,n\in\bN_0,m\neq n,\label{surf_harm:ortho}\\
  \intertext{and}
    \int D_n^\lambda(\eta\t x)D_n^\lambda(\xi\t x)\, \unif(\bS_d)(dx)& \ = \ \gamma_n^\lambda D_n^\lambda(\eta\t \xi),\quad \eta,\xi\in \Sd,\,n\in \bN_0; \label{surf_harm:convol} 
\end{align}
see e.g.~\citet[p.\,245]{ErdelyiII} and \citet[formula (1.14)]{SawJMVA}.

The mixing bases considered in what follows are of a very simple structure:
The densities of the base distributions are built with only two
special surface harmonics. To be precise, for $n\in\bN$ and real numbers $-1 \le \alpha \le +1$ 
the functions
\begin{equation}
  x\mapsto 1+\alpha D_n^\lambda(\eta\t x),\hskip2mm x\in \Sd,
\end{equation}
are $\unif(\Sd)$-densities of probability distributions $\Delta_{n,\eta,\alpha}^\lambda$ on $\Sd$. 
These distributions can be regarded as a multivariate generalization of the cardioid distributions
introduced by~\citet[p.\,302]{Jeffreys}; see also \citet[Section 3.5.5]{MardiaJupp}. Let $\Delta_{0,\eta,\alpha}^\lambda:=\unif(\bS_d)$.
Here we mainly deal with the special distributions 
$\Delta_{n,\eta}^\lambda:=\Delta_{n,\eta,1}^\lambda$, but see also Remark \ref{rem:extensions} (a) and
Proposition \ref{prop:selfmix} below.
So, for $n\in\bN$ the $\unif(\Sd)$-density of $\Delta_{n,\eta}^\lambda$ is simply the sum of the two special
surface harmonics  $H_{0,\eta}^\lambda\equiv 1$ and $H_{n,\eta}^\lambda$. 

With the mixing base $\Delta_{n,\eta}^\lambda,n\in\bN_0,$
given for each $\eta\in\Sd$, we obtain discrete mixture representations for all 
spherical distributions with densities of the form \eqref{eq:struct} that are not too far away from $\unif(\bS_d)$.
This may be seen as an 
instance of the perturbation approach discussed in the survey paper of~\citet{Pewsey} and,
indeed, the value of $\beta(\rho)$ may be interpreted as
a distance between $\nu_d$ and the measure with $\nu_d$-density $g_\rho$. 
By an \emph{ultraspherical mixing base} we mean a family 
$\{\Delta_{n,\eta}^\lambda:\, n\in\bN_0\}$.

The following general formula
is an immediate consequence of the above definitions and the expansion
in~\eqref{eq:expansion}.

\begin{proposition}\label{prop:gensphere}
Suppose that  $\beta_n(\rho)\ge 0$ for all $n\in\bN$ and
that $\beta(\rho)\le 1$. Let $\eta\in\bS_d$. Then $P_{\eta,\rho}$ has 
the discrete mixture expansion
\begin{equation}\label{eq:thmmix}
     P_{\eta,\rho}\, =\, \sum_{n=0}^\infty w_\rho(n) \, \Delta^\lambda_{n,\eta}, 
\end{equation}
with $w_\rho(0)=1-\beta(\rho)$ and
$w_\rho(n) =\beta_n(\rho)$ for all $n\in\bN$.
\end{proposition}  

Applying this construction to several specific families we have to take care of
the crucial condition $\beta(\rho)\le 1$, equivalently $w_\rho(0)\ge 0$. 
In each case, we obtain the mixing distribution and a range of $\rho$-values for the 
validity of the representation. Any distribution on $\bN_0$ may be written
as a mixture of unit mass at 0 and a distribution on $\bN$ and it turns out that 
the latter are occasionally from a standard family.
For a distribution on $\bN_0$ with mass function $w$ on $\bN_0$ such that $w(0)<1$ 
we call the distribution on $\bN$ with mass function  $n\mapsto w(n)/(1-w(0))$, $n\in\bN$,
its \emph{zero-truncated counterpart}. Some of the results stated in what follows
turn out to be simple consequences of Proposition \ref{prop:gensphere}. 

In the first theorem we consider two families of wrapped distributions, hence $d=1$ 
and $\lambda=0$. There are different notational conventions in the literature; 
here, we regard the wrapped distribution associated with a given distribution 
$\mu$ on (the Borel subsets of) the real line as the push-forward $\mu^T$ of $\mu$  
under the mapping $T:\bR\to\bS_1$, $x\mapsto (\cos(x),\sin(x))\t$. This is often
applied to location-scale families. Alternatively, the interval $[-\pi,\pi)$ is used instead of $\bS_1$
as the base set for the wrapped distribution. This means that with $X\sim \mu$ one deals with the
$[-\pi,\pi)$-valued random variable $X_0$ as the variable $X$ reduced modulo $2\pi$.   
If $X$ has the density $f$ with respect to the Lebesgue measure, then $X_0$ has the
$\unif\left ([-\pi,+\pi)\right )$-density 
$2\pi \sum_{n=-\infty}^{+\infty} f(s+2\pi n),\,s\in [-\pi,\pi)$. 
If the characteristic function $\varphi$ of $X$ is absolutely integrable, 
then $f$ is continuous and the Poisson summation formula applies, i.e.
\begin{equation}\label{PoisSum}
  2\pi \sum_{n=-\infty}^{+\infty} f(s+2\pi n)\,=\,\sum_{n=-\infty}^{+\infty} \varphi(n)\,e^{-ins},\quad s\in [-\pi,\pi);
\end{equation}    
see \citet[p.\,632]{Feller}. 
For example, the \emph{wrapped normal
distribution} $\WN_1(\eta,\rho)$ arises from the normal distribution $N(\alpha,\sigma^2)$ with mean
$\alpha$ and variance $\sigma^2$, where $\eta=(\cos(\alpha),\sin(\alpha))\t$ and
$\rho=\sigma^2$. Note that some authors use $\rho=2\sigma^2$, see, e.g. \citet{HartmanWatson}.   
With $X\sim N\bigl(\alpha,\sigma^2\bigr)$ we deduce from \eqref{PoisSum} that $X_0$ has the density
\begin{equation*}
1+2\sum_{n=1}^\infty e^{-n^2\rho/2} \cos n(s-\alpha),\quad s\in [-\pi,+\pi ),
\end{equation*}  
which means that $T\circ X\sim \WN_1(\eta,\rho)$ has the $\unif(\Seins)$-density
\begin{equation}\label{d:wnormal}
f_1^\WN(x|\eta,\rho)\;=\;1+2\sum_{n=1}^\infty e^{-n^2\rho/2}\,T_n(\eta\t x),\quad x\in \Seins .
\end{equation}
It is worthwhile to note that as $\rho\to 0$ the distribution $\WN_1(\eta,\rho)$ converges weakly to
the one-point mass distribution at $\eta$. This is in contrast to other distributions considered
here. For example, $\vMF_1(\eta,\rho),\, \Wat_1(\eta,\rho),\,\AG_1(\eta,\rho)$ all converge weakly to the uniform distribution $\unif(\Seins)$ as $\rho\to 0$. 
Also, as $\rho\to\infty$, the distribution $\WN_1(\eta,\rho)$ converges weakly to
$\unif(\Seins)$.

For the \emph{wrapped Cauchy distribution}
$\WC_1(\eta,\rho)$ we follow the definition given by~\citet{Pewsey}: If $X$ has a standard 
Cauchy distribution with density $x\mapsto 1/(\pi(1+x^2))$, $x\in\bR$, then we apply 
the wrapping procedure to $Y:=\sigma X+\alpha$, where $\sigma>0,\,\alpha\in\bR,$ and take
$\eta$ as in the wrapped normal case. For the scaling we use the parametrization
$\rho:=e^{-\sigma}\in (0,1)$ and augment this with the limiting uniform distribution at $\rho=0$.
Using \eqref{PoisSum} again we obtain that the distribution $\WC_1(\eta,\rho)$  
has the density
\begin{equation}\label{d:wcauchy}
  f_1^\WC(x|\eta,\rho)\;=\; 1+2\sum_{n=1}^\infty T_n(\eta\t x)\rho^n
  \;=\; \frac{1-\rho^2}{1-2\rho\,\eta\t x +\rho^2}, \quad x\in\bS_1;
\end{equation}
see also~\citet{Pewsey}. 

We write $\geo_{\bN_0}(n|p)=p(1-p)^n$, $n\in\bN_0$, for the probability mass function of
the geometric distribution on $\bN_0$ with parameter $p\in (0,1)$,
and $\geo_\bN$ for the mass function of its zero-truncated counterpart.
Recall that the function $\beta$ for a given distribution family is defined
in~\eqref{eq:def_beta}.

\begin{theorem}\label{thm:wrapped}
\emph{(a)} For the wrapped Cauchy distributions we have $\beta(\rho)= 2\rho/(1-\rho)$ and, 
for $0\le\rho\le 1/3$,
\begin{equation*}
   \WC_1(\eta,\rho) \, = \, (1-\beta(\rho))\, \unif(\bS_1)
               \, +\, \beta(\rho) \sum_{n=1}^\infty\geo_\bN(n|1-\rho)\,\Delta_{n,\eta}^0.
\end{equation*}

\vspace{.5mm}
\noindent
\emph{(b)} For the wrapped normal distributions we have 
$\beta(\rho)= 2\sum_{n=1}^\infty e^{-n^2 \rho/2}$. Let $\rho_0\approx 1.570818$ %0.785409$
be the unique solution of the equation $\beta(\rho)=1$. Then, with 
\begin{equation*}
   \br(n|\rho)\,  :=\, 2\beta(\rho)^{-1} e^{-n^2 \rho/2} %\exp\Bigl(-\frac{n^2 \rho}{2}\Bigr) 
      \quad\text{for all } n\in\bN,
\end{equation*} 
and $\rho\ge \rho_0$, it holds that
\begin{equation*}
   \WN_1(\eta,\rho) \, = \, (1-\beta(\rho))\, \unif(\bS_1)
               \, +\, \beta(\rho) \sum_{n=1}^\infty\br(n|\rho)\,\Delta_{n,\eta}^0.
\end{equation*}
\end{theorem}

We deal with the von Mises--Fisher families next. For these, we need variants of
the \emph{Skellam distribution} with parameter $\rho$. This distribution  arises as the 
distribution of $N_1-N_2$ where $N_1,N_2$ are independent random variables that 
both have the Poisson  distribution with parameter $\rho/2$; see~\citet{Irwin}, 
and see~\citet{Skellam} where the more general case with possibly different means for 
$N_1$ and $N_2$ is considered. 
In the one-dimensional case we need
the \emph{positive} Skellam distribution, %$\PSK(\rho)$
which is the conditional distribution of $|N_1-N_2|$ given that $N_1\neq N_2$. The 
associated probability mass function is given by
\begin{equation}\label{eq:pskrate}
     \psk(n|\rho)\, 
        =\, \frac{2e^{-\rho}I_n(\rho)}{1-e^{-\rho}I_0(\rho)}, \quad n\in\bN.
\end{equation}
For $d>1$ we use the \emph{generalized positive Skellam distribution}
with parameters $\kappa>0$ and $\tau>0$, with mass function
\begin{equation*}
     \gpsk(n|\kappa,\tau)\, :=\, 
           \Bigl( 1+\frac{n}{\kappa}\Bigr)\, \frac{(2\kappa)_n}{n!}
           \frac{2^\kappa\Gamma(\kappa+1)\tau^{-\kappa}e^{-\tau}I_{\kappa}(\tau)}
                {1-2^\kappa\Gamma(\kappa+1)\tau^{-\kappa}e^{-\tau}I_\kappa(\tau)}
           \frac{I_{\kappa+n}(\tau)}{I_\kappa(\tau)},\quad n\in\bN.
\end{equation*}  
It will turn out as part of the proof of the next result that this is indeed a 
probability mass function. We have
$\lim_{\kappa\to 0}\frac{1}{\kappa}\frac{(2\kappa)_n}{n!} = \frac{2}{n}$, which
gives $\lim_{\kappa\to 0}\gpsk(n|\kappa,\tau)=\psk(n|\tau)$ 
for all $n\in\bN$. Hence the
positive Skellam distribution appears as the limiting case  of the 
generalized positive Skellam distribution as $\kappa\to 0$.

\begin{theorem}\label{thm:sphericalMF} We consider the von Mises--Fisher families 
$\{\vMF_d(\eta,\rho):\, \rho\ge 0\}$, $\eta\in\bS_d$.

\noindent
\emph{(a)} 
If $d=1$ then $\beta(\rho)= e^\rho/I_0(\rho)-1$,  the
equation $\beta(\rho)=1$ has a unique finite positive solution
$\rho_0\approx 0.876842$, and for $0\le \rho\le \rho_0$ it holds that
\begin{equation*}
     \vMF_1(\eta,\rho)\; =\; (1-\beta(\rho))\, \unif(\bS_1)
                                          \, +\, \beta(\rho)\,
                                 \sum_{n=1}^\infty \psk(n|\rho)\,\Delta^0_{n,\eta}.
\end{equation*}
\noindent
\emph{(b)} If $d>1$ then
\begin{equation*}
    \beta(\rho) = \beta_\lambda(\rho) \, := \, \frac{\rho^\lambda e^\rho}
                             {2^\lambda \Gamma(\lambda + 1) I_\lambda(\rho)}
    \, - \, 1,
\end{equation*} 
the equation $\beta_\lambda(\rho)=1$ has a unique finite positive solution $\rho_0(\lambda)$, 
and for $0\le \rho\le \rho_0(\lambda)$ it holds that 
\begin{equation*}
        \vMF_d(\eta,\rho)
          \; =\; (1-\beta_\lambda(\rho))\,\unif(\bS_d)\, 
               +\, \beta_\lambda(\rho)
                          \sum_{n=1}^\infty \gpsk(n|\lambda,\rho)\Delta^\lambda_{n,\eta}.
\end{equation*} 
\end{theorem}

\begin{remark}\label{rem:extensions}
(a) The condition that $\beta_n(\rho) \ge 0$ for all $n\in\bN$ 
in Proposition~\ref{prop:gensphere} is satisfied 
in all of the above families, but it can easily be removed 
by an appropriate extension of the mixing base. 
For this,  let $\Delta^{\lambda,-}_{n,\eta}$ be the
distribution with density $x\mapsto 1-D_n^\lambda(\eta\t x)$, $x\in\bS_d$. Then 
the representation~\eqref{eq:thmmix} continues to hold if we take 
$Q_{n,\eta}=\Delta^{\lambda,-}_{n,\eta}$ and $w_\rho(n)=-\beta_n(\rho)$
whenever $\beta_n(\rho)<0$.

\vspace{1mm}
(b) The condition $\beta(\rho)\le 1$ in Proposition~\ref{prop:gensphere} holds if
$P_{\eta,\rho}$ is sufficiently close to the uniform distribution on
the sphere. If instead of $P_{\eta,\rho}$ we consider a mixture of this distribution
with the uniform, with enough weight on the latter, then the result is close enough
to the uniform,  and we again obtain a
mixture representation. In the von Mises--Fisher case with $d=1$,
for example, we get for $\rho>\rho_0$ and with 
$\alpha(\rho):= (1- 2e^{-\rho} I_0(\rho))/(1- e^{-\rho} I_0(\rho))$,
\begin{equation*}
   \alpha(\rho)\, \unif(\bS_1) \,+\,  (1-\alpha(\rho))\, \vMF_1(\eta,\rho)
                               \ = \ \sum_{n=1}^\infty \psk(n|\rho)\,\Delta^0_{n,\eta}.
\end{equation*}

(c) Is there a countable mixing base that represents \emph{all} isotropic spherical 
distributions with densities of the form $x \mapsto g(\eta\t x)$ with $\eta\in \Sd$   
and $g\in\bH_\lambda$? This may be rephrased in terms of the set of extremal points 
of a convex set in an infinite dimensional space. We refer to~\citet{BaGr6} for such
geometric aspects in general, and for the construction of tree-based mixing bases
that would lead to a positive answer for the set of all $g\in\bH_\lambda$ that are
Riemann integrable.
\end{remark}

The passage from an $L^2$-expansion~\eqref{eq:expansion} of $g_\rho$ to the 
mixture representation~\eqref{eq:thmmix} heavily relies on the nonnegativity of
the functions $1+D^\lambda_n$ (respectively $1-D^\lambda_n$ in part (a) above).
More generally, we may consider a mixing base $(Q_{n,\eta})_{n\in\bN}$ where 
the density of $Q_{n,\eta}$ is a polynomial of degree $n$ in $\eta\t x$.
This leads to the consideration of general linear combinations of 
ultraspherical polynomials; indeed, finding conditions for 
such polynomials to be nonnegative (on a given interval)
is an ongoing research topic, see~\citet{Askey}. 

We confine ourselves to an example with $\lambda=0$ and the Chebyshev polynomials.
A change of mixing base will obviously lead to a change
in the sequence of mixing coefficients. It turns out that this may lead to 
a representation of wider applicability. 

%In passing we note the resemblance to a similar phenomenon in the classical 
%summability theory of real sequences; see e.g.~[Hardy, Divergent Series].

\begin{example}\label{ex:lincombbase}
We define functions $g_\rho:[-1,1]\to \bR_+$, $0\le\rho < 1$, by
\begin{equation}\label{eq:gTuran}
   g_\rho(t) \, := \frac{(1-\rho t + \phi_\rho(t))^{1/2}}{2^{1/2}\phi_\rho(t)}, 
           \quad -1\le t\le 1,
\end{equation}
where $\phi_\rho(t):= (1-2\rho t +\rho^2)^{1/2}$. These can be written as
\begin{equation}\label{eq:gTuran2}
 g_\rho(t) = \sum_{n=0}^\infty\frac{(\nf{1}{2})_n}{n!}\, T_n(t) \,\rho^n,
\end{equation}
see \citet[p. 259]{MOS}. In particular, $\int_{-1}^1 g_\rho(t)\, dt = 1$, so that 
we may define a family of distributions $P_{\eta,\rho}$ on $\bS_1$ via their densities
$x\mapsto g_\rho(\eta\t x)$, $x\in\bS_1$.
  
We first derive an expansion in terms of the distributions 
$\Delta^0_{n\eta}$. With $t=1$ we get
\begin{align*}
    \beta(\rho)=\sum_{n=1}^\infty\frac{(\nf{1}{2})_n}{n!}\,\rho^n
                       \, = \, (1-\rho)^{-\nf{1}{2}}-1,
\end{align*}
and it follows that $\beta(\rho)\le 1$ if and only if $\rho\le 1-2^{-\nf{1}{2}}=:\rho_0$. 
We now introduce the \emph{zero-truncated negative binomial distribution} with parameters $r>0$, $p\in (0,1)$, and 
probability mass function
\begin{align*}
      \znb(n|r,p)=\frac{(r)_n}{n!}(1-p)^n\frac{p^r}{1-p^r},\quad n\in\bN.
\end{align*}  
Then \eqref{eq:gTuran2} leads to the discrete mixture representations
\begin{equation}\label{eq:firstrepr}
  P_{\eta,\rho}\, =\, \bigl(1-\beta(\rho)\bigr)\,\unif(\bS_1)\, +\, 
   \beta(\rho)\, \sum_{n=1}^\infty\znb(n|\nf{1}{2},1-\rho)\,\Delta^0_{n,\eta},
\end{equation}
for all $\eta\in\bS_1,\, 0< \rho \le \rho_0$.

On the other hand,~\citet{Turan} proved that, for all $n\in\bN_0$, 
\begin{equation*}
 \sum_{k=0}^n\frac{(\nf{1}{2})_k}{k!}\cos (k\vartheta) > 0\quad\text{for}~0<\vartheta<\pi.
\end{equation*}
In view of $T_k(\cos\vartheta)=\cos(k\vartheta)$ this can be used to obtain an 
alternative representation. For this, let $\Sigma_{n,\eta}$ be the distribution on $\bS_1$
with density $x\mapsto \sum_{k=0}^n \alpha_k T_k(\eta\t x)$, where
$\alpha_k:= (\nf{1}{2})_k/{k!}$ and $n\in\bN_0$. Then~\eqref{eq:gTuran2}, together with a summation
by parts, leads to 
\begin{equation}\label{eq:secondrepr}
  P_{\eta,\rho}\; = \;\sum_{n=0}^\infty\geo_{\bN_0}(n|1-\rho)\, 
                                    \Sigma_{n,\eta}
             \; =\;  (1-\rho)\,\unif(\bS_1) + \rho 
             \sum_{n=1}^\infty\geo_{\bN}(n|1-\rho)\,  \Sigma_{n,\eta},
\end{equation}
where $\geo_{\bN_0}$ and $\geo_\bN$ are as in Theorem~\ref{thm:wrapped}.
Both~\eqref{eq:firstrepr}
and~\eqref{eq:secondrepr} hold for all $\eta\in\bS_1$, but note that the range of
permissible $\rho$-values has increased from $0\le\rho\le 1-2^{-1/2}$ to the 
full interval $0\le \rho < 1$.
\end{example}

As pointed out earlier mixing families constructed with surface harmonics 
can be used with all distributions of the form~\eqref{eq:struct}
as long as these are sufficiently close to the uniform. The following result
gives a property 
which we interpret as \emph{self-mixing stability} of the distribution families
$\cD^\lambda_0:=\{\unif(\Sd)\}$, $\cD^\lambda_n:=\{\Delta^\lambda_{n,\eta,\alpha}:\, \eta\in\bS_d,\,-1\le \alpha\le +1\}$,
$n\in\bN$, and $\cD^\lambda:=\bigcup_{n\in\bN_0} \cD^\lambda_n$:
Mixing two elements of the same family results in a distribution that belongs
to this family as well. Additionally, mixing any two elements moves the 
mixing distribution closer to the 
uniform distribution. Generally, the mixing operation relates distributions with different 
location parameter $\eta\in\bS_d$ to each other.

\begin{proposition}\label{prop:selfmix}
Let $\gamma_n^\lambda$ be the constants defined in \eqref{eq:defgamma}.
  
\noindent
\emph{(a)} For all $n\in\bN_0$ and $\eta\in\bS_d$, $-1\le \alpha,\beta\le +1$,
\begin{align*}%\label{eq:selfmix1}
 \int \Delta^\lambda_{n,\zeta,\alpha}(A)\, \Delta^\lambda_{n,\eta,\beta}(d\zeta)
     &\; =\; \Delta^\lambda_{n,\eta,\gamma^\lambda_n\alpha \beta}(A)\\
     &\; =\;  (1-\gamma^\lambda_n)\, \unif(\bS_d)(A) 
       \, +\, \gamma^\lambda_n \, \Delta^\lambda_{n,\eta,\alpha \beta}(A)
\end{align*}
for all $A\in \cB(\bS_d)$.

\smallskip
\noindent
\emph{(b)} For all $n,m\in\bN_0$ with $n\not= m$, and all $\eta\in\bS_d$, $-1\le \alpha,\beta\le +1$,
\begin{equation*}%\label{eq:selfmix1}
 \int \Delta^\lambda_{m,\zeta,\alpha}(A)\, \Delta^\lambda_{n,\eta,\beta}(d\zeta)
     \; =\;   \unif(\bS_d)(A) 
                        \quad\text{for all } A\in\cB(\bS_d).
\end{equation*}
\end{proposition}

We recall that a \emph{probability kernel} from a measurable space $(E,\cE)$ to another 
measurable space $(F,\cF)$ is a function $Q:E\times \cF\to \bR$ that is 
$\cE$-measurable in its first and a probability measure on $(F,\cF)$ 
in its second argument. Given a probability measure $P$ on $(E,\cE)$ and a kernel 
$Q$ from $(E,\cE)$ to $(F,\cF)$ 
we define a probability measure $P\circ Q$ on $(F,\cF)$ by
\begin{equation}\label{eq:compo}
    P\circ Q(A) \, = \, \int Q(x,A)\, P(dx)\quad\text{for all } A\in\cF. 
\end{equation}
For a family $\{Q_x:\, x\in E\}$ of probability measures $Q_x$ on $(F,\cF)$ 
with the property that the map $x\mapsto Q_x(A)$ is
$\cE$-measurable for all $A\in \cF$, we may regard
$Q_{\bm \cdot}:E\times \cF \rightarrow\bR$ defined by 
$Q_{\bm \cdot}(x,A)=Q_x(A)$, $(x,A)\in E\times \cF$, as a Markov kernel
from~$(E,\cE)$ to $(F,\cF)$. Then~\eqref{eq:compo} reads
\begin{equation*}
  P\circ Q_{\bm\cdot}(A) \, = \, \int Q_x(A)\, P(dx)\quad\text{for all } A\in\cF, 
\end{equation*}
which may be interpreted as a mixing operation. 
For use in the next section we note that for a family
$\{P_\eta:\eta\in\Sd\}$ of probability measures $P_\eta$ on $(\Sd,\cB(\Sd))$
and a kernel $Q$ from~$(\Sd,\cB(\bS_d))$ to $(\bS_d,\cB(\bS_d))$ that are both isotropic 
in the sense defined by~\eqref{eq:isotropic} and~\eqref{ker:isotropic} respectively, 
the mixing results in an isotropic family again,
\begin{equation}\label{eq:isofamily}
   \bigl( P_{\eta}\circ Q \bigr)^U \, =\, P_{U\eta}\circ Q\quad \text{for all}~\eta\in\Sd~\text{and all}~U\in\OO(d+1).
\end{equation}
Further, for $\alpha=\beta=1$ the statements  
in Proposition~\ref{prop:selfmix} can be simply rephrased as
\begin{gather*}
    \Delta^\lambda_{n,\eta}\circ \Delta^\lambda_{n,\bm \cdot} 
          \; =\; (1-\gamma^\lambda_n)\, \unif(\bS_d) 
                         \, +\, \gamma^\lambda_n \, \Delta^\lambda_{n,\eta},\\
    \Delta^\lambda_{n,\eta}\circ \Delta^\lambda_{m,\bm \cdot}
          \; =\; \unif(\bS_d)\quad \text{if } n\not=m.
\end{gather*}
Using the self-mixing stability and the bilinearity of the operation defined
in~\eqref{eq:compo} we obtain a discrete mixture representation for the composition
of two families that both have a representation in terms of the ultraspherical mixing base.

\begin{proposition}\label{prop:compo}
Suppose that $P_\eta$ and $P'_\eta$, $\eta\in\bS_d$,  are distribution families on $\bS_d$
with discrete mixture representations $\, P_\eta=\sum_{n=0}^\infty w(n)\Delta^\lambda_{n,\eta}\,$
and $\,P'_\eta=\sum_{n=0}^\infty w'(n)\Delta^\lambda_{n,\eta}\,$ respectively.
Then
\begin{equation*}
     P_\eta\circ P'_{\bm \cdot}
     \; =\; \sum_{n=0}^\infty \tilde w(n) \, \Delta^\lambda_{n,\eta} ,
\end{equation*}
where $\; \tilde w(0):= 1 - \sum_{n=1}^\infty \tilde w(n)$ and 
$\; \tilde w(n):= \gamma^\lambda_n\, w(n)w'(n)\;$ for all $n\in\bN$.
\end{proposition}
 
This can be used to obtain the composition of wrapped Cauchy distributions. Indeed,
taken together, Theorem~\ref{thm:wrapped}\,(a) and Proposition~\ref{prop:compo} 
lead to
\begin{equation}\label{eq:WC-Levy}
  \WC_1(\eta,\rho)\circ \WC_1(\bm \cdot,\rho')\, = \, \WC_1(\eta,\rho\rho')
             \quad\text{for all }\rho,\rho'\le 1/3.
\end{equation}
It is worthwhile to point out that the restriction $\rho,\rho'\le 1/3$ in \eqref{eq:WC-Levy}
can be omitted. In fact, for all $0\le \rho,\rho'<1$, using
\eqref{eq:defgamma}, \eqref{surf_harm:ortho}, \eqref{surf_harm:convol}, and \eqref{d:wcauchy},
the density of
$\WC_1(\eta,\rho)\circ \WC_1(\bm \cdot,\rho')$ is easily calculated to be 
\begin{align*}
  \int f_1^\WC(x|\xi,\rho')\,f_1^\WC(\xi|\eta,\rho)\,\unif(\Seins)(d\xi)\ = \ f_1^\WC(x|\xi,\rho\rho'),\quad x\in\Seins.
\end{align*}  
In the next section this will be put into a wider context. 

\section{Discrete mixture representations and Markov processes}\label{sec:Markov}

Let $\{P_{\eta,\rho}:\,\eta\in\bS_d, \, \rho\in I\}$, $I\subset\bR_+$ an interval,
be a family  of distributions of the type considered in the previous sections. 
Below we briefly discuss
three different connections to Markov processes. First, for $\rho\in I$ fixed,
the corresponding subfamily may be regarded as a probability %(or Markov)  
kernel and thus induces
a Markov chain on spheres. Second, an isotropic  diffusion process
on $\bS_d$ leads to a family of the above type via its one-dimensional marginal
distributions, where $\eta$ and $\rho$ take over the role of starting point and
(transformed) time parameter respectively. Third, starting with a discrete mixture representation
we may find a discrete time Markov chain with marginal distributions equal to elements 
of the mixing base,
and thus obtain an almost sure representation of the family as the distributions
of the chain at random times.

\subsection{Random walks on spheres}\label{subsec:rw}
Let $\{Q_\eta:\, \eta\in\bS_d\}$ be a family of probability distributions that
leads to a Markov kernel as described at the end of Section~\ref{sec:ultra}.
Such kernels arise as transition probabilities of Markov processes. We may, 
for example, fix an $\eta\in\bS_d$ and define a Markov chain $(X_n)_{n\in\bN_0}$
with state space $(\bS_d,\cB(\bS_d))$ by the requirements
that $X_0\equiv\eta$ and that the distribution of $X_{n+1}$ conditionally on $X_n=\xi$
is given by~$Q_\xi$. For isotropic kernels each transition can be divided into
two steps that make use of the representation of $\bS_d$ by $[-1,1]\times \bS_{d-1}$
that also appeared in connection with~\eqref{eq:dens_hd} and~\eqref{eq:defC}. 
In the geometrically most familiar
case we consider the current position as the `north pole', then first choose a latitude 
and thereafter, independently, a longitude uniformly at random. The result is regarded 
as the new north pole. A generalization of this setup has been considered by~\citet{Bingham},
see also the references given there. 

The case $d=1$ is somewhat special as the wrapping procedure is a group homomorphism
from the additive group of real numbers into the multiplicative group $\bS_1$, regarded 
as a subset of $\bC$ and endowed with complex multiplication. Wrapping a random walk 
or a L\'evy process thus leads to processes with values in $\bS_1$ that have stationary
and independent increments, where the latter are now to be understood as ratios rather 
than differences. In fact, the location-scale family of Cauchy distributions arises as
the one-dimensional marginals of a specific L\'evy process, which 
gives~\eqref{eq:WC-Levy} after an appropriate rescaling of the variance parameter. 
A similar approach, now using Brownian motion on the real line, gives a corresponding 
statement for the wrapped normal distributions.

We collect some observations in the following result. Recall that the distribution
$\cL(X)$ of a Markov chain $X=(X_n)_{n\in\bN_0}$ with state space $(\bS_d, \cB(\bS_d))$
is a probability measure on the path space 
$(\bS_d^{\bN_0}, \cB\bigl( \bS_d^{\bN_0}\bigr))$ of the chain, where the $\sigma$-field
$\cB\bigl( \bS_d^{\bN_0}\bigr)$ is generated by the projections $\pi_k:\bS_d^{\bN_0}\to\bS_d$, 
$(x_n)_{n\in\bN_0}\mapsto x_k$, $k\in\bN_0$. Any measurable 
mapping $T:\bS_d\to \bS_d$ may be lifted to a mapping from and to paths by componentwise
application.

\begin{proposition}\label{prop:MCdiscrete}
Suppose that $Q$ is an isotropic kernel on $\bS_d$ and that $X^\eta=(X^\eta_n)_{n\in\bN_0}$ 
is a Markov chain with start at $\eta\in\bS_d$ and transition kernel $Q$. 

\vspace{0.5mm}
\noindent
\emph{(a)} The family $\{\cL(X^\eta):\, \eta\in\bS^d\}$ of probability measures on the
path space is isotropic in the sense that $\cL(X^\eta)^U=\cL(X^{U\eta})$ for all 
$\eta\in\bS_d$, $U\in\OO(d+1)$.

\vspace{0.5mm}
\noindent
\emph{(b)} The \emph{latitude process} $Y=(Y_n)_{n\in\bN_0}$, with $Y_n:=\eta\t X^\eta_n$
for all $n\in\bN_0$, is a Markov chain with state space $[-1,1]$ and start at $1$.

\vspace{0.5mm}
\noindent
\emph{(c)} Suppose that
$Q(\eta,\cdot\,) =\sum_{k=0}^\infty w(k)\Delta_{k,\eta}^\lambda$ for all $\eta\in \bS_d$.
Then the representation
\begin{equation}\label{eq:powerX}
 \cL(X_n^\eta)\, = \, \sum_{k=0}^\infty w_n(k)\Delta_{k,\eta}^\lambda
              \quad\text{for all }n\in \bN, \, \eta\in\bS_d,  
\end{equation}
holds with $w_1(k):=w(k)$ for all $k\in \bN_0$ and, for $n>1$,
\begin{equation*}
  w_n(k)\,:=\, (\gamma_k^\lambda)^{n-1}w(k)^n,\ k \in \bN,
            \quad w_n(0)\, := \, 1-\sum_{k=1}^\infty w_n(k).
\end{equation*}  
\end{proposition}

The fact that the geodesic distances from the starting point are again a Markov chain
is an instance of \emph{lumpability}, see~\citet{PitmanRogers} for a general discussion.
That the dependence on $\eta$ is lost in the lumping transition is part of the 
assertion of part~(b). Further, it follows from the cosine theorem for spherical triangles
that the transition kernel $Q_Y$ of $Y$ may be written as
\begin{equation}\label{eq:Bingham}
   Q_Y(y,\cdot\,) \, = \, \cL\bigl(yZ + U(1-y^2)^{1/2}(1-Z^2)^{1/2}\bigr)
\end{equation}
with $Z, U$ independent, $\cL(Z)=\cL(\eta\t X_1^\eta)$ and $\cL(U)=\nu_{d-1}$;
see also~\citet{Bingham}. Part~(c) shows that the mixing base introduced in 
Section~\ref{sec:ultra} is useful in the Markov chain context, and
\eqref{eq:powerX} may be seen as a discrete mixture representation of
the family $\{P_{\eta,n}:\, \eta\in\bS_d,\, n\in\bN\}$, with $P_{\eta,n}=\cL(X^\eta_n)$.
It implies that the marginal distributions of the associated latitude process are given by
\begin{equation*}
 \cL(Y_n)\, = \, \sum_{k=0}^\infty w_n(k)\, \mu_k^\lambda
              \quad\text{for all }n\in \bN, 
\end{equation*}
where $\mu_k^\lambda$ is the push-forward of $\Delta^\lambda_{k,\eta}$ under 
$x\mapsto \eta\t x$.

\subsection{Diffusion processes}\label{subsec:diff}
Let $X=(X_t)_{t\ge 0}$ be a homogeneous continuous time Markov process on the sphere 
with start at $\eta$, i.e.\ $\bP(X_0=\eta)=1$, and transition densities $p_t(x,y)$,
$t>0$, $x,y\in\bS_d$, that are isotropic in the sense that 
$p_t(Ux,Uy)=p_t(x,y)$ for all $t>0$, $x,y\in\bS_d$ and $U\in\OO(d+1)$.
Then the marginal distributions $\cL(X_t)$, $t\ge 0$, 
of the process may have a discrete mixture representation
of the type considered above, with $\rho$ related to time~$t$.

We sketch the basic argument, see also~\citet{Karlin}, and
then apply this in the context of the spherical Brownian motion $(B_t)_{t\ge 0}$.
For a general discussion of the latter we refer to~\citet[Section 7.15]{ItoMcKean} and~\citet[Example~3.3.2]{Hsu}. 
We note in passing that the 
marginal distributions of $X$ characterize the full distribution $\cL(X)$ of the process, 
in view of the Chapman--Kolmogorov equations and the invariance of the transition mechanism
under orthogonal transformations (clearly, $\OO(d+1)$ acts transitively on $\bS_d$).

Suppose that $X$ has transition densities $p_t(x,y)$ 
and that its infinitesimal generator $A$ has a discrete spectrum.
As transitions are isotropic it is enough to consider one specific starting value
$x=\eta$. For the Kolmogorov forward equations
\begin{equation*}
    \Bigl(\frac{\partial}{\partial t} p_{\bm.}(\eta,y)\Bigr)(t) 
        \, =\, \Bigl(Ap_t(\eta,\bm\cdot)\Bigr)(y)
\end{equation*}
we may try to find a family of basic solutions $\phi$ by a separation ansatz
$\phi(t,y)=f(t)g(y)$. This leads to 
\begin{equation*}
\frac{f'(t)}{f(t)}\, = \, \frac{(Ag)(y)}{g(y)}\, .
\end{equation*}
As the left and right hand side respectively depend on $t$ and $y$ only, we may hope that
\begin{equation*}
p_t(\eta,y) \, = \, \sum_{n=0}^\infty e^{-\omega_n t}\, \phi_{n,\eta}(y),
\end{equation*}
where $\omega_n$, $n\in\bN_0$, are the eigenvalues of the operator $A$, 
with eigenfunctions $\phi_{n,\eta}$. 

Recall that $\lambda=(d-1)/2$.

\begin{theorem}\label{thm:Brown}
Let $(B_t)_{t\ge 0}$ be the spherical Brownian 
motion on $\bS_d$, $d>1$, with start at $\eta\in \bS_d$. Let
\begin{equation*}
    \beta^\lambda(t) \, := \, \sum_{n=1}^\infty \Bigl(1+\frac{n}{\lambda}\Bigr)
                   \frac{(2\lambda)_n}{n!}\, e^{-n(n+2\lambda)t/2}
\end{equation*}
and let 
\begin{equation*}
          \br^\lambda_t(n):\, =\, \beta^\lambda(t)^{-1} 
                 \Bigl(1+\frac{n}{\lambda}\Bigr)
                   \frac{(2\lambda)_n}{n!}\, e^{-n(n+2\lambda)t/2}, \quad n\in\bN.
\end{equation*} 
Further, let $\,t^\lambda_0$ be the unique solution of the equation
$\beta^\lambda(t)=1$. Then, for $t\ge t^\lambda_0$,
\begin{equation*}
   \cL(B_t) \, = \, (1-\beta^\lambda(t))\, \unif(\bS_d)
         \, +\, \beta^\lambda(t) \sum_{n=1}^\infty\br_t^\lambda(n)\,\Delta_{n,\eta}^\lambda.
\end{equation*}
\end{theorem}

In view of the wrapping representation mentioned above the corresponding result for $d=1$ 
is contained in part~(b) of Theorem~\ref{thm:wrapped}. 

Let $B=(B_t)_{t\ge 0}$ be as in the theorem  and
let $Y=(Y_t)_{t\ge 0}$ with $Y_t=\eta\t B_t$, $t\ge 0$, 
be the associated latitude process. Then a polar decomposition shows that, 
for $t>0$ fixed, $B_t$ can be synthesized (in distribution) from $Y_t$ and an 
independent random variable $Z$ that is uniformly distributed on $\bS_{d-1}$. 
On the level of processes the conditional distribution of $B$ given $Y$ is a result 
known as the skew product decomposition of spherical Brownian motion; 
see~\citet[p.\,270]{ItoMcKean} and~\citet{Mija}.

A representation with a different mixing base, closer in spirit to the representations
in Section~\ref{sec:special}, has very recently been obtained by~\citet{Mija}.
The result is based on the authors' observation that for a spherical Brownian motion 
$(B^\lambda_t)_{t\ge 0}$ with start at $\eta\in\bS_d$, $d>1$, the rescaled latitude 
process $(Y_t)_{t\ge 0}$, 
$Y_t:=(1-\eta\t B^\lambda_t)/2$ for all $t\ge 0$, is a neutral Wright-Fisher diffusion
with both mutation parameters equal to $\lambda$. 
For these, a discrete mixture representation had earlier
been given by~\citet{JenkinsSpano}, see also~\citet{GriffithSpano}. Taken together, this leads to
\begin{equation}\label{eq:sphericalnormal}
   \cL(B^\lambda_t) \ = \  
           \sum_{n=0}^\infty w_t^\lambda(n)\,  \SBeta_d(\eta,1,n+1)\quad\text{for all } t>0,
\end{equation}
with
\begin{equation*} 
  w_t^\lambda(n) \ =\ \sum_{k=n}^\infty (-1)^{k-n}\frac{(d+2k-1)(d+n)_{k-1}}{n!(k-n)!}e^{-k(k+d-1)t/2}  
                  \quad\text{for all } n\in\bN_0,\ t>0, 
\end{equation*}
where the term $(d)_{-1}$ appearing with $n=k=0$ is defined to be $1/(d-1)$; see \citet[formula (5)]{JenkinsSpano}.  Interestingly,
the mixture coefficients turn out to be the individual probabilities associated with the 
marginal distributions of a particular pure death process $(Z_t)_{t>0}$, 
i.e.\ $w_t^\lambda(n)=P(Z_t=n)$. 
In contrast to our representation in Theorem~\ref{thm:Brown}
via surface harmonics no further restrictions on the time parameter are needed.
Moreover, there is also a fascinating
probabilistic interpretation, relating neutral Wright-Fisher diffusions to Kingman's
coalescent via moment duality; see~\citet{Mija} for the details.

\citet{MardiaJupp} call the $\cL(B^\lambda_t)$ Brownian motion distributions on $\Sd$;
\citet{Kent} regards $\cL(\eta\t B^\lambda_t)$ as a spherical normal distribution.
We adopt the notation of the latter. Remembering that $B^\lambda$
starts in $\eta\in\Sd$, we denote by $\SN_d(\eta,\rho)\,=\,\cL\bigl(B^\lambda_{1/\rho}\bigr)$
the spherical normal distribution
with parameters $\eta\in\Sd$ and $\rho>0$. Then, interestingly, by \eqref{eq:sphericalnormal}  
we have a
discrete mixture representation for the family $\{\,\SN_d(\eta,\rho):\eta\in\Sd, \rho>0\}\,$
with the same mixing base of spherical beta distributions $\SBeta_d(\eta,1,n+1)$, $n\in\bN_0$,
as obtained in Section~\ref{sec:general} for the von Mises--Fisher, the spherical Cauchy, 
and the angular Gaussian families. 

As explained at the beginning of this subsection, 
isotropic diffusion processes on the sphere may lead to a discrete
mixture representation for the family of their marginal distributions. 
Conversely, for a given family of the type considered in the previous sections, 
one might ask for a representation of its elements as the marginals 
of some diffusion process with values in $\bS_d$. The following result answers this 
question for the von Mises--Fisher distributions. 

\begin{theorem}\label{thm:vMF_Markov}
There is no homogeneous Markov process $X=(X_t)_{t\ge 0}$ on $\bS_d$ with the property that, 
for all $\eta\in\bS_d$,
\begin{equation}\label{eq:vMFfidis}
    \cL(X_t|X_0=\eta)\in\bigl\{\vMF_d(\eta,\rho):\, \rho > 0\bigr\} \quad\text{for all } t>0.
\end{equation}
\end{theorem}

Alternatively one might start with a diffusion on the ambient space $\bR^{d+1}$ and then use
the transition $x\mapsto x/\|x\|$ from $\bR^{d+1}$ to $\bS_d$. For example, if
$B=(B_t)_{t\ge 0}$ is a Brownian motion on $\bR^{d+1}$ with start at $\eta\in\bS_d$,
then $X=(X_t)_{t\ge 0}$, with $X_t:=\|B_t\|^{-1} B_t\,$ for all $t\ge 0$, 
represents the family $\{\AG_d(\eta,\rho):\,\rho>0\}$ in the sense that
\begin{equation*}
      \cL(X_t)=\AG_d\bigl(\eta,\rho(t)\bigr)\quad \text{ for all }t>0.
\end{equation*}
Here the bijection $\rho:(0,\infty)\to (0,\infty)$ is given by $\rho(t):=(2t)^{-1/2}$, $t>0$.

\begin{remark} 
We relate Theorem \ref{thm:vMF_Markov} to the infinite divisibility
statement for the von Mises--Fisher distributions obtained by \citet{Kent}. 
Interestingly, in both cases the proofs are based on the same series expansions \eqref{eq:seriesMF},
\eqref{eq:seriesMF2} of the densities $f_d^{\vMF}(\,\cdot\,|\eta,\rho)$ in terms of 
ultraspherical polynomials. An important ingredient of Kent's approach is the associative 
convolution algebra $(\mathscr F,\circ_d)$ on the space $\mathscr F$ of probability measures 
on $[-1,1]$ introduced by \citet{Bingham}, where the convolution $F_1\circ_d F_2$ of 
$F_1,F_2\in\mathscr F$ is defined to be the distribution of $S_1S_2 + \Lambda (1-S_1)^{\nf{1}{2}}(1-S_2)^{\nf{1}{2}}$, see also~\eqref{eq:Bingham}. 
Here $S_1,S_2,\Lambda $ are independent, 
$S_i\sim F_i$ for $i=1,2$, and $\Lambda\sim \nu_{d-1}$. For $d=0$ we take
$\nu_0$ to be the discrete uniform distribution on $\Snull:=\{-1,1\}$.    
With this definition
a probability measure $F\in \mathscr F$ is said to be \emph{$\circ_p$-infinitely divisible} 
if for each $m\in \bN$ there exists an
$F_m\in\mathscr F$ such that $F$ is equal to the $m$-fold convolution $F_m^{\circ_d m}$ of $F_m$.
As $F\in \mathscr F$ is uniquely determined by its
\emph{Fourier transform} $\varphi_F:\bN_0\rightarrow \bR$ defined by 
\begin{align*}
\varphi_F(m)\,=\,\int D_m^\lambda(t)\,F(dt)\quad\text{for}~m\in\bN_0, 
\end{align*}
and $\varphi_{F_1\circ_d F_2}=\varphi_{F_1}\varphi_{F_2}$ for all $F_1,F_2\in \mathscr F$,
it holds that $F$ is $\circ_d$-infinitely divisible if and only if for each $m\in\bN$ there
is a Fourier transform $\varphi_{F_m}:\bN_0\rightarrow \bR$ of an
$F_m\in\mathscr F $ such that $\varphi_F=\varphi_{F_m}^m$ for all $m\in\bN$.
The distribution $\vMF_d^*(\rho)$ of $\eta\t X$, with $X\sim \vMF_d(\eta,\rho)$,
has the $\nu_d$-density $\frac{\rho^d}{2^d\Gamma(\lambda+1)I_\lambda(\rho)}\exp(\rho t),\,t\in [-1,1]$.
\citet{Kent} showed that $\vMF_d^*(\rho)$ is $\circ_d$-infinitely divisible. 
In fact, Kent even gives an interesting representation of the distributions $F_m$ such that $\vMF_d^*(\rho)=F_m^{\circ_d m}$ based on the spherical Brownian motion $(B_t)_{t\ge 0}$ on 
$\bS_d$ with start at $\eta\in \bS_d$: He showed that
for each $m\in\bN$ there exists an absolutely continuous distribution $G_m$ on the positive half-line
such that $\eta\t B_{T_m}\sim F_m$, where $T_m\sim G_m$ is independent of $(B_t)_{t\ge 0}$. 
Consequently,
\begin{equation}\label{eq:convolution}
\vMF_d^*(\rho)\,=\,\cL\left (\eta\t B_{T_m}\right )^{\circ_d m}\quad \text{for all}~m\in\bN. 
\end{equation}
In order to lift this from the unit interval to the sphere let $\eta\in\Sd$ be fixed 
and let $\mathscr P_\eta$ be the family of distributions $P$ on $\Sd$ that
are axially symmetric with respect to $\eta$, i.e.\ $P^U=P$ for each $U\in \OO(d+1)$ 
with $U\eta=\eta$, and $P^U$ the push-forward of $P$ under $U$. In order to carry over
the convolution operation $\circ_d$ to $\mathscr P_\eta$, a spherical addition $\oplus$ 
on $\Sd$ is defined.
Assign to each $x\in \Sd$ an element $U_x\in \OO(d+1)$ in such a way that $U_x\eta=x$ for all $x\in\Sd$ and
that $x\rightarrow U_x$ is a measurable injection (a measurable embedding) from $\Sd$ into $\OO(d+1)$.
Then, for $x,y\in \Sd$, let $x\oplus y:=U_xy$. For $P_1,P_2\in \mathscr P_\eta$, with independent 
$\Sd$-valued random vectors $X_i\sim P_i$, $i=1,2$, the $\star_d$-convolution $P_1\star_d P_2$ of $P_1$ and $P_2$
is defined to be the distribution of $X_1\oplus X_2$. Kent showed that $P_1\star_d P_2\in \mathscr P_\eta$, and
with $F_i$ as the distribution of $\eta\t X_i$, $i=1,2$, the cosine theorem for spherical 
triangles leads to $\eta\t (X_1\oplus X_2)\sim F_1\circ_d F_2$, i.e.
\begin{equation*}
 \cL\left (\eta\t (X_1\oplus X_2)\right )\,=\, \cL(\eta\t X_1) \circ_d \cL(\eta\t X_2).
\end{equation*} 
From \eqref{eq:convolution} it now follows that 
$\vMF_d(\eta,\rho)\,=\,\cL\bigl(B_{T_m}\bigr)^{\star_d m}$ for all $m\in\bN$. 
\end{remark}

\subsection{Almost sure representations}\label{subsec:asrepr}
The basic relation~\eqref{eq:repr} connects $\{P_\theta:\, \theta\in\Theta\}$
to the distributions $W_\rho$ on $\bN_0$ and $Q_{n,\eta}$ on $\bS_d$. 
Isotropy means that we may consider the $\eta$-part of the parameter
as fixed. In this situation, if $(N_\rho)_{\rho\ge 0}$ and $(X_n)_{n\in\bN_0}$ are 
independent stochastic processes such that
\begin{equation*}
   \cL(N_\rho)=W_\rho\text{ for all }\rho\ge 0,\quad 
   \cL(X_n)=Q_{n,\eta}\text{ for all }n\in\bN_0,
\end{equation*}
then
\begin{equation}\label{eq:asrepr}
     P_{\eta,\rho} =\cL(X_{N_\rho})\quad\text{for all }\rho\ge 0.
\end{equation}
Equation~\eqref{eq:asrepr} may be regarded as an almost sure representation 
of the distributional equation~\eqref{eq:repr}. Classical examples of such 
almost sure representations  are the Skorohod coupling in connection with
distributional convergence, see e.g.~\citet[Theorem 3.30]{Kall}, and the representation 
of a sequence of uniformly distributed permutations on the sets $\{0,1\ldots,n\}$, 
$n\in\bN$, by the Chinese restaurant process, see e.g.~\citet[Section 3.1]{Pitman}. 
Of course, such representations are most useful if the successive variables
are close to each other (and not simply chosen to be independent). 
Our aim here are representations of the type~\eqref{eq:asrepr} for the distribution families 
considered above. For this, we first formalize the polar decomposition. 

Let $\eta\in\bS_d$, we assume that $d>1$. Recall from~\eqref{eq:defC} that 
$C_d(\eta,y)=\{z\in\bS_d:\, \eta\t z = y\}$ and let
\begin{equation*}
  n_\eta:\bS_d\setminus\{-\eta,\eta\}\to C_d(\eta,0),\quad 
           x \mapsto \frac{x-(\eta\t x)\eta}{\|x-(\eta\t x)\eta\|}, 
\end{equation*} 
be the normalized projection onto the orthogonal complement of the linear 
subspace of $\bR^{d+1}$ spanned by $\eta$. This can be extended to the whole of 
the sphere by choosing some arbitrary of $\xi\in\bS_d$ as the value of $n_\eta(\pm \eta)$. 
Then an inverse of the polar decomposition $x\mapsto (\eta\t x,n_\eta(x))$ is given by
\begin{equation*}
       \Psi_\eta :[-1,1]\times C_d(\eta,0)\to\bS_d,\quad   
                      (y,z) \mapsto  y\eta+ (1-y^2)^{1/2} z,
\end{equation*}
in the sense that $\; \Psi_\eta(\eta\t x,n_\eta(x)) = x$ for all $x\in\bS_d$, and
a random vector $X$ with values in $\bS_d$ may be written as
\begin{equation}\label{eq:tangent-normal1}
 X =\Psi_\eta\bigl(\eta\t X,n_{\eta}(X)\bigr).
\end{equation}
Let $Q_\eta$ be a distribution on $(\bS_d,\cB(\bS_d))$ that has a density $f$
with respect to $\unif(\bS_d)$ which can be written as $f_\eta(x)=g(\eta\t x)$, $x\in\bS_d$,
see \eqref{eq:basic}. If $X\sim Q_\eta$ then 
$\eta\t X$ and $n_\eta(X)$ are independent, $\eta\t X$ has the distribution $\nu_{d;g}$
with $\nu_d$-density $g$, and $n_\eta(X)\sim \unif\left (C_d(\eta,0)\right )$; see
\citet[p.\,92]{Watson} and \citet[p.\,169]{MardiaJupp}. So, conversely, if we have independent random variables 
$Y\sim \nu_{d;g}$ and $Z\sim \unif\left (C_d(\eta,0)\right )$, then
\begin{equation}\label{eq:tangent-normal2}
  X\,=_{\mathcal D}\,\Psi_\eta (Y,Z)\,=\,Y\eta\,+\,(1-Y^2)^{1/2}Z.
\end{equation}  
\citet[p.\,161,\,p.\,169]{MardiaJupp} call \eqref{eq:tangent-normal1}
and the distributional version \eqref{eq:tangent-normal2} the 
\emph{tangent-normal decomposition}. In the past this decomposition has been successfully applied by many authors treating different problems in directional statistics,
see e.g. \citet{Saw,SawBiom,SawJMVA}, \citet{GPV2020}, and \citet{Ulrich}.
In practice, a random variable $X$ with distribution $Q_\eta$ is simply obtained as follows.
Suppose first that $\eta$ is equal to $e_1=(1,0,\ldots,0)\t$, the first unit vector in the canonical basis of $\bR^{d+1}$.
With $e_1$ as `north pole' the polar representation takes on a particularly simple
form: From $y\in[-1,1]$ and $z=(z_1,\ldots,z_d)\t\in\bS_{d-1}$ we get $x\in\bS_d$
by $x=\Phi(y,z)$ with 
\begin{equation}\label{eq:def_polar}
\Phi(y,z)\,:=\, \bigl(y\, ,\,(1-y^2)^{1/2}\, z_1,\ldots,(1-y^2)^{1/2}\,z_d\bigr)\t\,.
\end{equation}
Then, starting with a random variable $Y\sim \nu_{d;g}$ and another random variable $Z\sim\unif(\bS_{d-1})$
independent of $Y$, we obtain an $\bS_d$-valued random variable $X$ with distribution $Q_{e_1}$
via $X:=\Phi(Y,Z)$.
For a general $\eta\in\bS_d$ we use that $Q_\eta$ is the push-forward $Q_{e_1}^U$ of $Q_{e_1}$
under the mapping $x\mapsto Ux$, where $U\in\OO(d+1)$ is such that $\eta=Ue_1$,
i.e.~$U$ has $\eta$ as its first column. Then, defining $\Phi_\eta(y,z)=U\Phi(y,z)$ 
it follows that $X:=\Phi_\eta(Y,Z)\sim Q_\eta$; see also \citet{Saw} and \citet{Ulrich} for this construction.

Starting with a random variable $Y$ that is almost surely equal to 1, it follows that
the random variable $X=\Phi_\eta(Y,Z)$ is almost surely equal to $\eta$.
This means that from a Markov chain 
$(Y_n)_{n\in\bN_0}$ with state space $[-1,1]$ starting in 1, where for $n\in\bN$ the distribution of $Y_n$ 
is the push-forward of $Q_{n,\eta}$ under $x\mapsto \eta\t x$, and a \emph{single} 
random variable $Z\sim\unif(\bS_{d-1})$, we obtain a Markov chain $(X_n)_{n\in\bN_0}$
with the desired one-dimensional marginal distributions by putting $X_n:=\Phi_\eta(Y_n,Z)$ 
for all $n\in\bN_0$. Clearly, $(Y_n)_{n\in\bN_0}$ is then the latitude process associated with $(X_n)_{n\in \bN_0}$.

This reduces the first step in an almost sure construction~\eqref{eq:asrepr} 
to finding a Markov chain
on $[-1,1]$ with prescribed marginals. For the second step we require a suitable 
integer-valued process $N=(N_\rho)_{\rho\ge 0}$ with marginal distributions $W_\rho$. 
One general possibility is the quantile transformation,
which can also be used to construct a Skorohod coupling for real random variables:
With $U\sim\unif(0,1)$ we obtain a random variable $X$ with distribution function $F$
via $X:=F^{-1}(U)$, where $F^{-1}(u):=\inf\{t\in\bR;\, F(t)\ge u\}$ for $0<u<1$. 
If $F=F_\rho$ is the distribution function associated with $W_\rho$  the paths of 
the process $N=(N_\rho)_{\rho\ge 0}$ constructed in this way depend on the relations
between the distribution functions for different $\rho$'s. In particular, if the
distributions $W_\rho$ 
are stochastically monotone, meaning that
\begin{equation}\label{eq:stoch:monotone}
  1-F_\rho(x) \le 1- F_{\rho'}(x)\quad \text{for all}~x\in\bR
\end{equation}  
whenever $\rho\le \rho'$, then the paths of $N$ are increasing.
It is well known that this stochastic monotonicity applies to arbitrary distributions $W_\rho,W_{\rho'}$
with monotone likelihood ratio. To be precise, defining a likelihood ratio of $W_{\rho'}$ with respect to
$W_\rho$ to be a $\cB(\bR)$-measurable function $L_{\rho,\rho'}:\bR\rightarrow [0,\infty ] $ such that
$W_{\rho}(L_{\rho,\rho'}<\infty)\,= 1\; $ and 
\begin{align*}
  W_{\rho'}(A)& \ = \ \int_A L_{\rho,\rho'}(x)\, W_{\rho}(dx) +  W_{\rho'}(\{L_{\rho,\rho'}=\infty\}\cap A)\hspace{3mm}\text{for all}~A\in \cB(\bR),
\end{align*}  
the distributions $W_\rho,W_{\rho'}$ with $\rho<\rho'$ have monotone likelihood ratio if there exists an 
increasing function $h_{\rho,\rho'}:\bR\rightarrow [0,\infty ]$ such that
\begin{equation*}
  L_{\rho,\rho'}\,=\,h_{\rho,\rho'}~(W_\rho+W_{\rho'})\text{-almost everywhere}.
\end{equation*}
Note that if $f_{\rho},f_{\rho'}$ are densities of $W_\rho,W_{\rho'}$ with respect to some
$\sigma$-finite measure, then 
\begin{equation*}
  L_{\rho,\rho'}^*\,=\frac{f_\rho'}{f_\rho}\,1(f_\rho>0) \ + \ \infty \,1(f_\rho=0,f_{\rho'}>0)
\end{equation*}
is a special version of the likelihood ratio of $W_{\rho'}$ with respect to $W_\rho$; here
$1(\cdot)$ denotes the indicator function. In fact, \eqref{eq:stoch:monotone} holds if
$W_\rho,W_{\rho'}$ with $\rho<\rho'$ have monotone likelihood ratio; see, e.g. \citet[Satz 2.28]{Witting}. 

Again, we consider the special case of the von Mises--Fisher 
distributions in some detail.

\begin{example}\label{ex:vMFasrepr}
Let $\eta\in\Sd$. In order to translate Theorem~\ref{thm:mixrepr_vMF},
see also Example~\ref{ex:bsp}\,(c), into an almost sure representation for the family 
$\{\vMF_d(\eta,\rho):\,\rho\ge 0\}$ we need 
random variables $Y_n$ with $Y_n\sim\Beta_{[-1,1]}(d/2,d/2 + n)$ for all $n\in\bN_0$. 
To this end let $V$, $W_i$, $i\in\bN_0$, be independent random variables with $V\sim\Gamma(d/2,1)$, $W_0\sim\Gamma(d/2,1)$, and $W_i\sim\Exp(1)$ for all $i\in\bN$. 
Then, using the well known connection between beta and gamma distributions, 
see e.g.~\citet{JohnsonKotz}, we obtain independent random variables $B_0\sim\Beta(d/2,d/2)$,
$B_n\sim \Beta(d+n-1,1)$, $n\in\bN$, via
\begin{equation*}
     B_0:= \,  \frac{V}{V+W_0}, \quad 
            B_n:=\frac{V+W_0+\ldots+W_{n-1}}{V+W_0+ \cdots + W_n}\ \ \text{ for all } n\in\bN,
\end{equation*}
and products 
\begin{equation*}
    \tilde Y_n:= \prod_{i=0}^n B_i\; = \; \frac{V}{V+W_0+\ldots + W_n}\; 
           \sim\; \Beta(d/2,d/2+n)\quad\text{for all }n\in\bN_0.
\end{equation*}
The transformation $Y_n:=1-2 \tilde Y_n$, $n\in\bN_0$, now gives the desired sequence $Y=(Y_n)_{n\in\bN_0}$.
Moreover, $\tilde Y_{n+1}= \tilde Y_n B_{n+1}$ implies that
\begin{equation}\label{eq:Markov} 
  Y_{n+1} = 1 - (1-Y_n)\, B_{n+1}.
\end{equation} 
As $B_{n+1}$ is independent of $Y_n$ this shows that $Y$ is a Markov chain.

Suppose now that $(N_\rho)_{\rho\ge 0}$ is a stochastic process
with $N_\rho\sim \CHS(d/2,d,2\rho)$ for all $\rho\ge 0$. Let $Z\sim\unif(\bS_{d-1})$
be independent of the variables $V$ and $W_i$, $i\in\bN$. Then $(X_{\rho})_{\rho\ge 0}$
with 
\begin{equation}\label{eq:vMFas2}
   X_\rho\,:=\, \Phi_\eta\bigl(Y_{N_\rho} , Z\bigr) 
                     \quad\text{for all }\rho\ge 0,
\end{equation}
and $\Phi_\eta$ as in the remarks preceding the example,
has the desired property that $X_\rho\sim\vMF_d(\eta,\rho)$ for all $\rho\ge 0$.

For the construction of the counting process we use the quantile transformation.
The likelihood ratios turn out to be
\begin{equation*}
         \frac{\chs(n|d/2,d,2\rho')}{\chs(n|d/2,d,2\rho)} \, =\, 
         \frac{{}_1F_1(d/2;d;\rho)}{{}_1F_1(d/2;d;\rho')}
                 \, \frac{(2\rho')^n}{(2\rho)^n}
\end{equation*}
which, as a function of $n$, is increasing whenever $\rho < \rho'$.  
As explained above, this shows that the paths of $(N_\rho)_{\rho\ge 0}$ are increasing. 
From \eqref{eq:Markov} it follows that $Y_n\ge Y_{n-1}$ for all $n\in\bN$.  
Taken together we see that we have found an almost sure representation \eqref{eq:vMFas2} with
a process $(Y_{N_\rho})_{\rho\ge 0}$ that has increasing paths. 
\end{example}

Some comments are in order. Obviously, almost sure representations are generally
not unique. In the first step of Example~\ref{ex:vMFasrepr}, we could use a sequence 
$(Y_n)_{n\in\bN_0}$ of independent random variables with  $Y_n\sim \Beta_{[-1,1]}(d/2,d/2 + n)$ for all $n\in\bN_0$, 
or we could use the quantile transformation to obtain suitable variables $Y_n$ as functions
of one single $U\sim\unif(0,1)$ (in fact, the
corresponding likelihood ratios would be increasing). Similar to~\eqref{eq:chiquadrat} 
in the classical case,  the representation~\eqref{eq:vMFas2} strikes a
structural middle ground in this spectrum from no dependence at all to total 
dependence between the variables of interest. Also, the Markov chain featuring
in the denominator of 
\begin{equation*}
  Y_{N_\rho}=1-\frac{2V}{V+\sum_{n=0}^{N_\rho}W_n}
\end{equation*}
has some resemblance to the sum appearing
in~\eqref{eq:chiquadrat}. In~\citet[Remark 4\,(a)]{BaGr5} we found a 
discrete mixture representation for the non-central family of hyperbolic secant 
distributions that may similarly written as a function of a Markov chain of
this type.
%Further, we should mention that other results in Section~\ref{sec:special} 
%can similarly be used to obtain almost sure representations for other families
%of distributions on spheres.

Returning to the general situation, we may regard the right hand side
of~\eqref{eq:chiquadrat} or~\eqref{eq:asrepr} as a representation of a continuous time
stochastic process $Z= (Z_t)_{t\ge 0}$ by independent processes $X=(X_n)_{n\in\bN_0}$
and $N=(N_t)_{t\ge 0}$ via $Z_t=X_{N_t}$ for all $t\ge 0$. As already mentioned in the introduction, equality of the marginal distributions is considerably weaker than 
equality of the distributions of the processes. For example, the representation 
in Example~\ref{ex:vMFasrepr} leads to a process that moves by jumps from $\eta$ on 
a fixed great circle through $\eta$ towards the equator in a piecewise constant 
manner. Loosely speaking, a discrete mixture representation on the process level is
only possible for processes of the pure jump type. However, if the time parameter of 
the base process is $X$ continuous too then we obtain a 
connection to a famous group of results, known as \emph{skew product decompositions}. 
For example, with $(X_t)_{t\ge 0}$ a Brownian motion on $\bR^{d+1}$ starting at 
$a\not=0$ and $R_t:=\|X_t\|$, $t\ge 0$, we have 
$R_t^{-1} X_t = B_{N_t}$ for all $t\ge 0$,
where $(B_t)_{t\ge 0}$ is a Brownian motion on $\bS_d$ starting
at $\eta:=\|a\|^{-1} a$, $N_t$ is given implicitly
by $\int_0^{N_t} R_s^{-2}\, ds=t$, and $B$ and
$N$ are independent. In particular, this leads to a representation of the distribution
of $X_t/\|X_t\|$, $t\ge 0$, which is a family of spherical distributions,  
as a \emph{continuous} mixture. In contrast to discrete mixture representations these
seem to be less suitable for simulation.

\section{Proofs}\label{sec:proofs}

\subsection{Proof of Theorem~\ref{thm:mixrepr_vMF}}%\label{subsec:proof_vMF}
(a) We first simplify the norming constant in~\eqref{eq:normSBeta} for the $d$-dimen\-sional 
spherical beta distribution 
with parameters $p=1$ and $q=n+1$, $n\in\bN_0$. Using the duplication formula
\begin{equation}\label{eq:dupl}
    \Gamma(\nf{1}{2})\Gamma(2z)\,=\, 2^{2z-1}\Gamma(z)\Gamma(z+\nf{1}{2}), \quad z>0,
\end{equation}
for the gamma function we obtain
\begin{align}
  c_d(1,n+1)\; &=\; 2^{-(n+2\lambda)}\frac{\Gamma(\nf{1}{2})\Gamma(n+2\lambda+1)}
                                           {\Gamma(\lambda+1)\Gamma(n+\lambda+\nf{1}{2})}\notag \\
            &=\; 2^{-(n+2\lambda-1)}\frac{(2\lambda+1)_n}{(\lambda+\nf{1}{2})_n}
                                    \frac{\Gamma(\nf{1}{2})\Gamma(2\lambda)}
                                               {\Gamma(\lambda)\Gamma(\lambda+\nf{1}{2})}\label{eq:betanorming}\\
            &=\; 2^{-n}\frac{(2\lambda+1)_n}{(\lambda+\nf{1}{2})_n},\notag
\end{align}
and it follows that the density of $\SBeta_d(\eta,1,n+1)$ can be written as
\begin{equation}\label{eq:betadist_1}
  f_d^{\SBeta}(x|\eta,1,n+1)\; 
             =\; 2^{-n}\frac{(2\lambda+1)_n}{(\lambda+\nf{1}{2})_n}
                        \left (1+\eta\t x\right )^{n},\quad x\in \Sd.
\end{equation}
We now write
\begin{equation*}
      \exp\bigl(\rho\,\eta\t x\bigr) \; 
           = \; \exp(-\rho)\exp\bigl(\rho\, (1+\eta\t x)\bigr)
\end{equation*} 
and use the expansion
\begin{align*}
    \exp\bigl(\rho\, (1+\eta\t x)\bigr)\ 
         &=\ \sum_{n=0}^\infty\frac{\rho^n}{n!}\bigl(1+\eta\t x\bigr)^n\\
    &=\sum_{n=0}^\infty\frac{(2\rho)^n}{n!} \frac{(\lambda+\nf{1}{2})_n}{(2\lambda+1)_n}\,  
                                    f_d^{\SBeta}(x|\eta,1,n+1)(x),\quad x\in \bS_d,
\end{align*}
together with the identity
\begin{equation*}
     e^{-\rho}{}_{1}F_1(\lambda+\nf{1}{2};2\lambda+1;2\rho)\; 
              =\; 2^\lambda\Gamma(\lambda+1)I_\lambda(\rho)/\rho^\lambda
            \end{equation*}
to obtain 
\begin{equation*}
     f_d^{\vMF}(x|\eta,\rho)\; 
             =\; \sum_{n=0}^\infty \chs(n|\lambda+\nf{1}{2},2\lambda+1,2\rho)\;
                                          f_d^{\SBeta}(x|\eta,1,n+1),\quad x\in \bS_d.
\end{equation*}

In order to prove uniqueness suppose that
\begin{equation*}%\label{eq:uniqueMF}
  \sum_{n=0}^\infty v_n\, \SBeta_d(\eta,1,n+1)\; 
              = \; \sum_{n=0}^\infty w_n\, \SBeta_d(\eta,1,n+1)
\end{equation*}
for two sequences $v=(v_n)_{n\in\bN_0},w=(w_n)_{n\in\bN_0}$ of non-negative real numbers with
sum~$1$. Passing to the respective push-forwards under $x\mapsto \eta\t x$ this leads to 
the equality of the (continuous) densities,
\begin{equation*}
 \sum_{n=0}^\infty v_n\, \frac{n+1}{2^{n+1}}\, (1+y)^n\; 
              = \; \sum_{n=0}^\infty w_n\, \frac{n+1}{2^{n+1}}\, (1+y)^n,\quad -1<y<1,
\end{equation*}
hence the sequences $v$ and $w$ are equal to each other.

\vspace{0.5mm}
(b) From $d=2\lambda+1$, $1-\frac{4\rho}{(1+\rho)^2} = \frac{(1-\rho)^2}{(1+\rho)^2}$, and \eqref{eq:betanorming} it follows that
\begin{align*}
  f_d^\CII(x|\eta,\rho)&\,=\, \left (\frac{1-\rho^2}{(1+\rho)^2}\right )^{2\lambda+1}
  \left (1-\frac{2\rho}{(1+\rho)^2}\left (1+\eta\t x\right )\right )^{-(2\lambda+1)}\\
  &\,=\, \left (\frac{1-\rho}{1+\rho}\right )^{2\lambda+1}
               \sum_{n=0}^\infty \frac{(2\lambda + 1)_n}{n!}
                    \left (\frac{2\rho}{(1+\rho)^2}\right )^n
                         \left (1+\eta\t x\right )^n\\
  &\,=\, \left (\frac{(1-\rho)^2}{(1+\rho)^2}\right )^{\lambda+\nf{1}{2}}
          \sum_{n=0}^\infty \frac{(2\lambda + 1)_n}{n!}
            \left (\frac{4\rho}{(1+\rho)^2}\right )^n\,
            \frac{(\lambda+1/2)_n}{(2\lambda+1)_n}\;f_d^{\SBeta}(x|\eta,1,n+1)\\
    &\,=\,\sum_{n=0}^\infty \nb\left (n|\lambda+\nf{1}{2},(1-\rho)^2/(1+\rho)^2\right )\;f_d^{\SBeta}(x|\eta,1,n+1),\quad x\in \bS_d. 
\end{align*}
The proof of the uniqueness is similar to that of part (a).   

\vspace{0.5mm}
(c) Noticing $(d+1)/2=\lambda +1$, and
\begin{equation*}
  {}_2F_1\left (\lambda+\nf{1}{2},\lambda+1;2\lambda+1;4\rho/(1+\rho)^2\right )\,=\,\frac{(1+\rho)^{2\lambda+1}}{1-\rho},
\end{equation*}
see \citet[p.\,39]{MOS}, we get arguing as in part (b) that
\begin{align*}
  f_d^\CII(x|\eta,\rho)&\,=\,\frac{1-\rho^2}{(1+\rho)^{2(\lambda+1)}}\left (1-\frac{2\rho}{(1+\rho)^2}\left (1+\eta\t x\right )\right )^{-(\lambda+1)}\\
  &\,=\,\frac{1-\rho}{(1+\rho)^{2\lambda+1}}\sum_{n=0}^\infty \frac{(\lambda +1)_n}{n!}
                         \left (\frac{4\rho}{(1+\rho)^2}\right )^n\,\frac{(\lambda + 1/2)_n}{(2\lambda +1)_n}\,
     \;f_d^{\SBeta}(x|\eta,1,n+1) \\
  &\,=\,\sum_{n=0}^\infty \hs\left (n|\lambda+\nf{1}{2},\lambda+1,2\lambda+1,4\rho/(1+\rho)^2\right )\;f_d^{\SBeta}(x|\eta,1,n+1),\quad x\in \Sd.  
\end{align*}

\subsection{Proof of Theorem~\ref{thm:mixrepr_Wat}}\label{subsec:proof_Wat}
(a) Straightforward manipulations give
\begin{align*}
  f_d^{\Wat}(x|\eta,\rho)\ &=\ \frac{1}{{}_1F_1(\nf{1}{2};\lambda+1;\rho)}
                     \; \sum_{n=0}^\infty\frac{\rho^n (\eta\t x)^{2n}}{n!}\\
             &=\ \sum_{n=0}^\infty \chs(n|\nf{1}{2},\lambda+1,\rho)\, f_d^{\SP}(x|\eta,n).
\end{align*}

\vspace{0.5mm}
(b) As at the beginning of the proof of Theorem~\ref{thm:mixrepr_vMF}, with $p=q=n+1$, $n\in\bN_0$, 
the general expression~\eqref{eq:normSBeta} for the norming constants can be simplified to
\begin{equation*}
  c_d(n+1,n+1)\, =\, \frac{(\lambda+1)_n}{(\lambda+\nf{1}{2})_n}.
\end{equation*}
Hence the density for the associated spherical beta distribution may be written as  
\begin{equation}\label{betadist:2}
  f_d^{\SBeta}(x|\eta,n+1,n+1)\; =\; \frac{(\lambda+1)_n}{(\lambda+\nf{1}{2})_n}
        \bigl(1-(\eta^t x)^2\bigr)^n,\quad  x\in \Sd.
\end{equation}
For $\rho<0$ and $\eta\in \Sd$ we have
\begin{equation*}
    e^{\rho \,(\eta\t x)^2}
       \ =\ e^\rho e^{-\rho (1-(\eta\t x)^2)}
    \ =\ e^\rho\sum_{n=0}^\infty \frac{(-\rho)^n}{n!}
              \frac{(\lambda+\nf{1}{2})_n}{(\lambda+1)_n}
                         \, f_d^{\SBeta}(x|\eta,n+1,n+1)
\end{equation*}
for all $x\in\Sd$. Because of
\begin{equation*}
      e^{-\rho}{}_{1}F_1(\nf{1}{2};\lambda+1;\rho)\, 
                        =\, {}_{1}F_1(\lambda+\nf{1}{2};\lambda+1;-\rho)
\end{equation*}
(see, e.g. \citet[p.\,267]{MOS}) it follows that
\begin{align*}
    f_d^{\Wat}(x|\eta,\rho)\; =\; \sum_{n=0}^\infty \chs(n|\lambda+\nf{1}{2},\lambda+1,-\rho)          \, f_d^{\SBeta}(x|\eta,n+1,n+1)
\end{align*}
for all $\eta\in \Sd$, $\rho<0$. 

In both cases, it is easy to adapt the uniqueness argument from the von Mises--Fisher context
to the Watson situation.

\subsection{Proof of Lemma~\ref{lem:dpc}}
We write $\Gamma(\delta,\alpha)$ for the gamma distribution with shape parameter $\delta>0$,
scale parameter $\alpha>0$ and Lebesgue density 
$t\mapsto \Gamma(\delta)^{-1}\alpha^\delta t^{\delta-1}\exp(-\alpha t)$, $t>0$.
Let $V$ be a positive random variable with $V\sim\Gamma(\delta,\nf{1}{2})$. 
Then, for $k\in \bN_0$,
\begin{align*}
  E\left (V^{k/2}\exp\bigl(-2^{\nf{1}{2}}\,\tau\, V^{\nf{1}{2}}\bigr) \right ) \, = \, \frac{\Gamma(k+2\delta)}{2^{\delta -1}\Gamma(\delta)}
  \,e^{\tau^2/2}\,D_{-(k+2\delta)}(2^{\nf{1}{2}}\tau).
\end{align*}
Writing
\begin{displaymath}
  \frac{(\delta-\nf{1}{2})_k}{(2\delta-1)_k}=\frac{\Gamma(\delta-\nf{1}{2}+k)}{\Gamma(\delta -\nf{1}{2})}\frac{\Gamma(2\delta-1)}{\Gamma(2\delta -1 +k)} 
  \end{displaymath}
and using the duplication formula~\eqref{eq:dupl} we obtain
\begin{equation}\label{eq:rel2gamma}
  \dpc(k|\delta,\tau) =  \frac{1}{k!}(2^{\nf{1}{2}}\tau)^k 2^k\frac{(\delta-\nf{1}{2})_k}{(2\delta-1)_k}\,e^{-\tau^2}
       E\left (V^{k/2}\exp\bigl(-2^{\nf{1}{2}}\,\tau\, V^{\nf{1}{2}}\bigr)\right ).
\end{equation}
Thus,
\begin{align*}
  \sum_{k=0}^\infty \dpc(k|\delta,\tau)
        & \ = \ e^{-\tau^2}\sum_{k=0}^\infty  \frac{1}{k!}(2^{\nf{1}{2}}\tau)^k
                    2^k\frac{(\delta-\nf{1}{2})_k}{(2\delta-1)_k}
                     E\left (V^{k/2}\exp\bigl(-2^{\nf{1}{2}}\,\tau V^{1/2}\bigr) \right )\\
        & \ = \ e^{-\tau^2} E\left (\exp\bigl(-2^{\nf{1}{2}}\,\tau V^{1/2}\bigr)
           \,{}_{1}F_1(\delta-\nf{1}{2};2\delta-1;2^{\nf{3}{2}}\tau V^{1/2})\right ).
\end{align*}
Using
\begin{align*}
  e^{-z}\,{}_{1}F_1(\delta-\nf{1}{2};2\delta-1;2z) \ = \ I_{\delta-1}(z)\Gamma(\delta)(z/2)^{-(\delta-1)},
\end{align*}
for $z\in\bR$, see e.g. \citet[p. 265, formula (10)]{ErdelyiI}, and
\begin{displaymath}
  I_{\delta-1}(z)\Gamma(\delta)(z/2)^{-(\delta-1)} \ = \ \sum_{k=0}^\infty \frac{(z/2)^{2k}}{k!\Gamma(k+\delta)}
\end{displaymath}
we obtain
\begin{align*}
  \sum_{k=0}^\infty \dpc(k|\delta,\tau)& \ = \ \Gamma(\delta)e^{-\tau^2}
                                         \sum_{k=0}^\infty \frac{(\tau^2/2)^k}{k!\Gamma(k+\delta)}E(V^k),
\end{align*}
which in view of  $E(V^k)=\frac{\Gamma(k+\delta)}{\Gamma(\delta)}2^k$ gives 
$\sum_{k=0}^\infty \dpc(k|\tau,\delta)=1$.

\subsection{Proof of Theorem~\ref{thm:mixrepr_Saw}}%
Using \eqref{Saw:rep0} with $S\sim\chi^2_{d+1}$ and writing 
\begin{equation*}
\exp\left ( 2^{\nf{1}{2}}\rho S^{\nf{1}{2}}\,\eta\t x\right ) \ 
    = \ \exp\left ( 2^{\nf{1}{2}}\rho S^{\nf{1}{2}}\,(1+ \eta\t x)\right )
         \exp\left (-2^{\nf{1}{2}}\rho S^{\nf{1}{2}}\right )
\end{equation*}
we see that
\begin{align*}
  f_d^{\AG}(x|\eta,\rho)
          & = \sum_{n=0}^\infty (1+\eta\t x)^n\frac{\left (2^{\nf{1}{2}} 
            \rho\right )^n}{n!} e^{-\rho^2} E\left (S^{n/2}\exp(-2^{\nf{1}{2}}
                                                     \rho S^{\nf{1}{2}})\right ) \\
          & = \sum_{n=0}^\infty f_d^{\SBeta}(x|\rho,1,n+1)\,
                 \frac{\left(2^{\nf{1}{2}}\rho\right)^n}{n!}\,
                           2^n\frac{(\lambda+1/2)_n}{(2\lambda+1)_n}\,e^{-\rho^2} 
                     E\left (S^{n/2}\exp(-2^{\nf{1}{2}}\rho S^{\nf{1}{2}})\right ).
\end{align*}
As $\chi_{d+1}^2=\Gamma(\lambda+1,\nf{1}{2})$,  the asserted discrete mixture representation
now follows with~\eqref{eq:rel2gamma}. The uniqueness of the representation is obtained as in
the proofs of Theorems~\ref{thm:mixrepr_vMF} and~\ref{thm:mixrepr_Wat}.

\subsection{Proof of Theorem~\ref{thm:wrapped}}

(a) We use Proposition~\ref{prop:gensphere}. We have
$\beta(\rho)=\sum_{n=1}^\infty2\rho^n=\frac{2\rho}{1-\rho}\le 1$ if and only if
$\rho\le \nf{1}{3}$. By \eqref{d:wcauchy},
the density $f_1^\WC(\cdot|\eta,\rho)$ of the wrapped Cauchy distribution $\WC_1(\eta,\rho)$
can be written as 
\begin{align*}
     f_1^\WC(x|\eta,\rho)\; =\; 1-\beta(\rho)\, +\,2\sum_{n=1}^\infty \bigl (1+T_n(\eta\t x)\bigr )\rho^n,\quad x\in\bS_1.
\end{align*}
The representation now follows easily.

\vspace{1mm}
(b) Using \eqref{d:wnormal} we can write  the density $f_1^\WN(\cdot|\eta,\rho)$ 
of the wrapped normal distribution $\WN_1(\eta,\rho)$ as 
\begin{align*}
     f_1^\WN(x|\eta,\rho)\ =\ 1-\beta(\rho)\, +\, 
                      2\,\sum_{n=1}^\infty e^{-n^2\rho/2}\,\bigl(1+T_n(\eta\t x)\bigr),
                 \quad x\in\bS_1.
\end{align*}
The function $\beta$ is easily seen to be continuous and strictly decreasing, with unique 
solution $\rho_0\approx 1.570818$ of the equation $\beta(\rho)=1$. 
From Proposition~\ref{prop:gensphere} we thus obtain a 
discrete mixture expansion for the family 
$\{\WN_1(\eta,\rho):\, \eta\in\bS_1, \rho\ge \rho_0\}$ with
mixing base elements $\Delta^0_{n,\eta}$ and 
weights $w_\rho(n)= 2e^{-n^2\rho/2}$, $n\in\bN$.

\subsection{Proof of Theorem~\ref{thm:sphericalMF}}
(a) The density $f_1^\vMF(\cdot|\eta,\rho)$ of the von Mises--Fisher distribution
$\vMF_1(\eta,\rho)$ can be written as
\begin{equation}\label{eq:seriesMF}
     f_1^\vMF(x|\eta,\rho)\; =\; \frac{\exp(\rho\eta\t x)}{I_0(\rho)}
          \;=\; 1\,+\,2\,\sum_{n=1}^\infty \frac{I_n(\rho)}{I_0(\rho)}\,T_n(\eta\t x),
                \quad x\in\bS_1, 
\end{equation}
see e.g.~\citet[formula 9.6.34]{AbramSteg}. In fact, the representation \eqref{eq:seriesMF} corresponds to
the Fourier series expansion of the density of the unique random angle in $[-\pi,\pi)$
associated with $X\sim \vMF_1(\eta,\rho)$. To be precise, let for simplicity $\eta = (1,0)\t$ and let
$\Theta\in [-\pi,\pi)$ be the unique random angle such that $X=(\cos \Theta,\sin \Theta)\t$. Then the density
$h_\rho(\theta)=\frac{\exp(\rho \cos \theta)}{2\pi I_0(\rho)},\,\theta\in [-\pi,\pi),$ of $\Theta$ with respect to
the Lebesgue measure on $[-\pi,\pi)$ has the absolutely convergent Fourier series expansion
\begin{equation*}
     h_\rho(\theta)\; 
            =\; 1\,+\,2\,\sum_{n=1}^\infty \frac{I_n(\rho)}{I_0(\rho)}\,\cos(n\theta),
                \quad \theta \in [-\pi,\pi); 
\end{equation*}
see, \citet[formula (1.1a)]{Kent}.
In particular, 
$\beta(\rho) = 2\sum_{n=1}^\infty \frac{I_n(\rho)}{I_0(\rho)}
    =e^\rho/I_0(\rho)-1$.
We have  $\lim_{\rho\to 0}\beta(\rho)=0,\,\lim_{\rho\to\infty}\beta(\rho)=\infty$,
see \citet[formula 9.7.1]{AbramSteg}, and
\begin{equation*}
\frac{d}{d\rho}\beta(\rho)=I_0(\rho)^{-2}e^\rho\left (I_0(\rho)-I_1(\rho)\right )\,>\,0
                \quad \text{for all }\rho>0,
\end{equation*}  
see \citet{Soni}, hence the function $\beta$ is strictly increasing. Taken together 
this implies that the equation $\beta(\rho)=1$ has a unique finite positive root $\rho=\rho_0\approx 0.876842$. Proposition~\ref{prop:gensphere} now leads to
the discrete mixture representation  
\begin{align*}
     \vMF_1(\eta,\rho)\ &=\
         \left (2-e^\rho/I_0(\rho)\right )\unif(\bS_1)\, +\, 
                                       \left (e^\rho/I_0(\rho)-1 \right )
                              \sum_{n=1}^\infty \frac{2e^{-\rho}I_n(\rho)}
                                        {1-e^{-\rho}I_0(\rho)}\Delta^0_{n,\eta}\\
                        &=\ (1-\beta(\rho))\, \unif(\bS_1)
                                          \, +\, \beta(\rho)\,
                                 \sum_{n=1}^\infty \psk(n|\rho)\,\Delta^0_{n,\eta}\,.
\end{align*}

%\vspace{.8mm}
(b) 
For $\kappa>0,\,\tau>0$, and $-1\le t\le 1$, we have 
\begin{equation}\label{Gegenbauer:exp}
           e^{\tau t} \; =\;  2^\kappa \Gamma(\kappa)\tau^{-\kappa}\sum_{n=0}^\infty
                                  (\kappa+n)I_{\kappa + n}(\tau)C_n^\kappa(t),
\end{equation}
see \citet[p.\,227]{MOS} or \citet[Section 7]{Kent}. Hence the density $f_d^\vMF(\cdot|\eta,\rho)$ of 
$\vMF_d(\eta,\rho)$ can be written as
\begin{equation}\label{eq:seriesMF2}
  f_d^\vMF(x|\eta,\rho)\;=\; 1+\sum_{n=1}^\infty 
                        \Bigl(1+\frac{n}{\lambda}\Bigr) \frac{(2\lambda)_n}{n!}
           \frac{I_{\lambda+n}(\rho)}{I_\lambda(\rho)}D_n^\lambda(\eta\t x),
                         \quad x\in\Sd. 
\end{equation}
With $t=1$ in \eqref{Gegenbauer:exp} we get, after some algebra, 
\begin{equation*}
     1\; =\; \sum_{n=1}^\infty \Bigl(1+\frac{n}{\kappa}\Bigr) \frac{(2\kappa)_n}{n!}
               \frac{2^\kappa\Gamma(\kappa+1)\tau^{-\kappa}e^{-\tau}I_{\kappa}(\tau)}
                   {1-2^\kappa\Gamma(\kappa+1)\tau^{-\kappa}e^{-\tau}I_\kappa(\tau)}
              \frac{I_{\kappa+n}(\tau)}{I_\kappa(\tau)}.
\end{equation*}
This implies that $\gpsk(\cdot|\kappa,\tau)$ is a probability mass function.
Further, 
\begin{equation*} 
  \beta_\lambda(\rho) \ = \  \frac{\rho^\lambda e^\rho}{2^\lambda \Gamma(\lambda +1 ) I_\lambda(\rho)}\;-\; 1
                     \  = \ \sum_{n=1}^\infty \Bigl(1+\frac{n}{\lambda}\Bigr)
                      \frac{I_{\lambda +n}(\rho)}{I_\lambda(\rho)}C_n^\lambda(1).
\end{equation*}
Arguing as in part (a) we find that the equation
$\beta_\lambda(\rho)=1$ has a unique finite positive root $\rho=\rho_0(\lambda)$. 
For the family $\{\vMF_d(\eta,\rho):\,\eta\in\Sd,0<\rho\le \rho_0(\lambda)\}$
this finally gives the discrete mixture representation
\begin{equation*}
        \vMF_d(\eta,\rho)
          \; =\; (1- \beta_\lambda(\rho))\, \unif(\bS_d)\, 
               +\, \beta_\lambda(\rho)\,
                          \sum_{n=1}^\infty \gpsk(n|\lambda,\rho)\Delta^\lambda_{n,\eta}.
\end{equation*} 

\subsection{Proof of Proposition~\ref{prop:selfmix}}
The statements are a consequence of the basic formulas \eqref{surf_harm:ortho}
and \eqref{surf_harm:convol}.

\subsection{Proof of Proposition~\ref{prop:compo}}
All terms involved are positive, so there are no convergence issues, and 
the representation follows from the bilinearity of the mixture operation.

\subsection{Proof of Proposition~\ref{prop:MCdiscrete}}
(a) The distribution of $X$ is determined by the distributions of the vectors
$(X_0,\ldots,X_n)$, $n\in\bN_0$. We may thus use induction, and~\eqref{eq:isofamily} 
provides the necessary argument for the induction step.

\vspace{0.5mm} 
(b) Given a north pole $\eta$ we obtain a partitioning of $\bS_d$ into the sets
$C_\eta(y)=\{x\in\bS_d:\ \eta\t x= y\}$ with the same latitude $y\in [-1,1]$. 
Isotropy implies that the transition mechanism interacts with the function that
maps the points of $\bS_d$ to their latitude in the manner required by Dynkin's
criterion, see \citet[Theorem 10.13]{Dynkin}. As the action of $\OO(d+1)$ on $\bS_d$
is doubly transitive, for unit vectors $\eta_1,\eta_2\in \Sd$
with the property that $\eta\t \eta_1=\eta\t \eta_2$ there exists some $U\in\OO(d+1)$
such that $U\eta = \eta$ and $U\eta_1 = \eta_2$. The isotropy of $Q$ then implies that
for all $A\in \cB([-1,1])$ 
\begin{align*}
  Q(\eta_2,\{x\in\Sd:\eta\t x\in A\})&= Q(U\eta_1,\{x\in\Sd:\eta\t x\in A\})\\
                                     &= Q(\eta_1,\{U\t x\in\Sd:\eta\t x\in A\})\\
                                     &= Q(\eta_1,\{y\in\Sd:\eta\t Uy\in A\})\\
                                     &= Q(\eta_1,\{y\in\Sd:\eta\t y\in A\}).
\end{align*}

\vspace{0.5mm} 
(c) This follows easily on using induction and Proposition~\ref{prop:compo}.

\subsection{Proof of Theorem~\ref{thm:Brown}}

The transition probabilities for the spherical Brownian motion in dimension~$d\ge 2$ are given 
in~\citet{HartmanWatson},
\begin{align*}
       p_t(x,y)\; =\; 1+\sum_{n=1}^\infty \Bigl(1+\frac{n}{\lambda}\Bigr)
                  \frac{(2\lambda)_n}{n!}e^{-n(n+2\lambda)t/2}\,
                              D_n^\lambda(x\t y),\quad x,y\in\Sd,\, t>0,
\end{align*}
see also \citet{Karlin}. Let $P_t(\cdot,x)$ be the distribution on $\bS_d$ with density
$y\rightarrow p_t(x,y)$, let
\begin{align*}
     \beta_\lambda(t)\,:=\, \sum_{n=1}^\infty\Bigl(1+\frac{n}{\lambda}\Bigr)
           \frac{(2\lambda)_n}{n!}e^{-n(n+2\lambda) t/2}\quad\text{for }t>0,
\end{align*}
and let $t_\lambda>0$ be the unique positive real number such that
 $\beta_\lambda(t_\lambda)=1$. Then, with $\br_d(\cdot|t)$ as in the theorem,
we get the discrete mixture representations
\begin{align*}
         P_t(\cdot,\eta)\, =\, \left (1-\beta_\lambda(t)\right )\unif(\bS_d)\, 
                             +\, \beta_\lambda(t)\sum_{n=1}^\infty
                                            \br_d(n|t)\Delta^\lambda_{n,\eta}.
\end{align*}

\subsection{Proof of Theorem~\ref{thm:vMF_Markov}}
Suppose that, on the contrary, $(X_t)_{t\ge 0}$ is a homogeneous Markov process 
with the property that there exists a function $\kappa:(0,\infty)\rightarrow (0,\infty)$ 
such that for all $\xi\in\Sd$ it holds that if $\bP(X_0=\xi)=1$ then
\begin{equation*}
  \cL(X_t)=\vMF_d\left (\xi,\kappa(t)\right )\quad\text{for all}~t>0.
\end{equation*}  
Because of the Markov property we would then have
positive parameters $\rho,\rho',\rho''$ such that
\begin{equation}\label{eq:MFcompo}
 \vMF_d(\eta,\rho) \circ \vMF_d(\bm \cdot,\rho') = \vMF_d(\eta,\rho'').
\end{equation}
To see that this cannot be true, we note that the density of 
$\vMF_d(\eta,\rho) \circ \vMF_d(\bm \cdot,\rho')$ is given by
\begin{align*}
  h(x)\ := \ \int f_d^\vMF(x|\xi,\rho')\,f_d^\vMF(\xi|\eta,\rho)\,\unif(\Sd)(d\xi),\quad x\in\bR^d.
\end{align*}  
We now refer to \citet[Section 7]{Kent} and the proof of Theorem \ref{thm:sphericalMF}, where it is shown that the density
$f_d^\vMF(\bm\cdot|\eta,\rho)$ 
can be written as 
\begin{align*}
  f_d^\vMF(x|\eta,\rho) \ = \ \sum_{n=0}^\infty\frac{1}{\gamma_n^\lambda}\frac{I_{\lambda+n}(\rho)}{I_\lambda(\rho)}\,D_n^\lambda(\eta\t x)\quad\text{for all }x\in\Sd.
\end{align*}
Expressing $f_d^\vMF(\bm\cdot|\xi,\rho')$ correspondingly, and using
\eqref{surf_harm:ortho} and \eqref{surf_harm:convol}, we obtain
\begin{align*}
  h(x)\ = \ \sum_{n=0}^\infty\frac{1}{\gamma_n^\lambda}\frac{I_{\lambda+n}(\rho)}{I_\lambda(\rho)}
  \frac{I_{\lambda+n}(\rho')}{I_\lambda(\rho')} \,D_n^\lambda(\eta\t x)
         \quad\text{for all }x\in\Sd.
\end{align*}
Hence by \eqref{eq:MFcompo} we would have 
\begin{equation}\label{eq:contradict}
  \frac{I_{\lambda+n}(\rho)}{I_\lambda(\rho)}\frac{I_{\lambda+n}(\rho')}{I_\lambda(\rho')}
\ = \ \frac{I_{\lambda+n}(\rho'')}{I_\lambda(\rho'')}\quad\text{for all}~n\in\bN_0.
\end{equation}
From \eqref{BesselI} it is easily seen that for each $\tau>0$ the asymptotic relation
\begin{align*}
  I_{\lambda+n}(\tau)\sim 2^{-(\lambda+n)}\Gamma(\lambda+n+1)^{-1} \tau^{\lambda+n}\quad\text{as}~n\to\infty
\end{align*}
holds. Thus, 
\begin{align*}
  \log I_{\lambda+n}(\tau) = -(\lambda+n)\log 2 - \log \Gamma(\lambda+n+1) + (\lambda+n)\log \tau + o(1)\quad\text{as}~n\to\infty.
\end{align*}
By Stirling's formula, 
\begin{align*}
  \log \Gamma(\lambda+n+1) = (\lambda+n+1/2)\log\,n-n+\frac{1}{2}\log\,(2\pi)+o(1)\quad\text{as}~n\to\infty.
\end{align*}
From this we deduce that
\begin{equation*}
  \lim_{n\to\infty}\frac{\log \frac{I_{\lambda+n}(\tau)}{I_\lambda(\tau)}}{n\log n}= -1,
\end{equation*}  
which is in contradiction to~\eqref{eq:contradict}.

\vspace{1cm}
\noindent
{\bf \Large Acknowledgment}

\bigskip
\noindent
The authors would like to thank the reviewers for helpful suggestions.

\vspace{1cm}
\noindent
{\bf \Large Statements and Declarations}

\bigskip
\noindent
The authors declare that they have no conflict of interest. 

\vspace{0.5cm}
%\bibliographystyle{spbasic} 

%\vspace{3cm}

\end{document}